\documentclass[12pt]{amsart}
\usepackage[utf8]{inputenc}
\usepackage[margin=2.1cm]{geometry}
\usepackage{amsmath,amsfonts,accents}	
\usepackage{amsthm}
\usepackage{graphicx}
\usepackage{cancel}
\usepackage{multicol}
\usepackage{subfig}
\usepackage{multirow}
\usepackage[usenames,dvipsnames]{color}
\usepackage{charter}
\usepackage{color}
\usepackage{comment}
\usepackage{hyperref}

\usepackage{threeparttable}

\usepackage{algorithm}
\usepackage{algpseudocode}
\makeatletter
\def\BState{\State\hskip-\ALG@thistlm}
\makeatother

\newenvironment{rcases}
  {\left.\begin{aligned}}
  {\end{aligned}\right\rbrace}

\newtheorem{theorem}{Theorem}
\newtheorem{lemma}[theorem]{Lemma}
\newtheorem{proposition}[theorem]{Proposition}
\newtheorem{definition}[theorem]{Definition}
\title[Contour Parametrization via Mean Curvature Flows]{Contour Parametrization  via Anisotropic Mean Curvature Flows}

\author{P. Su\'arez-Serrato}
\thanks{Support for this project came from DGAPA-UNAM PAPIIT  IN102716 and UC-MexUS CN-16-43 grants, and the PASPA program from DGAPA-UNAM}

\address{Department of Mathematics, University of California, Santa Barbara, {\it on leave from}, Instituto de Matem\'aticas, Instituto de Matem\'aticas, {\sc Universidad Nacional Aut\'onoma de M\'exico, Mexico City} }

\email{pablo@im.unam.mx}

\author{ E.I. Vel\'azquez Richards}

\address{Instituto de Matem\'aticas, {\sc Universidad Nacional Aut\'onoma de M\'exico, Mexico City}}

\date{\today}

\begin{document}

\begin{abstract}
 We present a new implementation of anisotropic mean curvature flow for contour recognition. Our procedure couples the mean curvature flow of planar closed smooth curves, with an external field from a potential of point-wise charges. This coupling constrains the motion when the curve matches a picture placed as background. We include a stability criteria for our numerical approximation.
\end{abstract}

\maketitle

\tableofcontents

\section*{Introduction}
Many applied and theoretical problems have been studied through the analysis of manifold deformations~\cite{chopp1993,chow,GerhardPolden1999}. The description of these deformations by an evolution equation imposed via a geometric quantity are referred to as {\it geometric flows}. Applications include, for example, the growth of crystals, the modeling of fluids, and digital image recognition. Since their conception there has been continued interest in the development of numerical approximations to these flows~\cite{osher-level}.

The difficulty in analyzing these flows numerically depends on the geometric quantity in evolution (e.g. curvature, metric tensor, the manifold itself). In particular, the mean curvature flow (MCF) deforms a hypersurface in the normal direction $\mathbf{\hat{n}}$ with a speed proportional to its mean curvature $H$. This flow has an associated quasilinear parabolic equation in terms of an immersion $X$ of the hypersurface into the ambient manifold:
\begin{equation}\label{eq:MCF}
\frac{\partial X}{\partial t}=-H\mathbf{\hat{n}}
\end{equation}

The numerical methods employed to solve this equation are classified as Eulerian or Lagrangian methods, depending on the discrete representation of the surface or curve in evolution. In Eulerian methods the evolution is tracked by values at fixed positions on a gridded ambient. Conversely, in Lagrangian methods, the object in evolution is tracked explicitly through the position of its points. In consequence, Lagrangian methods have two advantages: they require smaller data storage than Eulerian methods, and the solution is computed explicitly; although generally their error estimates are difficult to estimate precisely. Eulerian methods are now a powerful and frequently used technique for mean curvature flow applications since the development of the Level Set Method~\cite{osher-level,sethian99}. To our knowledge, there are very few examples of Lagrangian methods that approach this problem~\cite{sevcovic2011, kimura}, our work adds to this list.

\begin{figure} [H]
  \centering
  \includegraphics[scale=0.3]{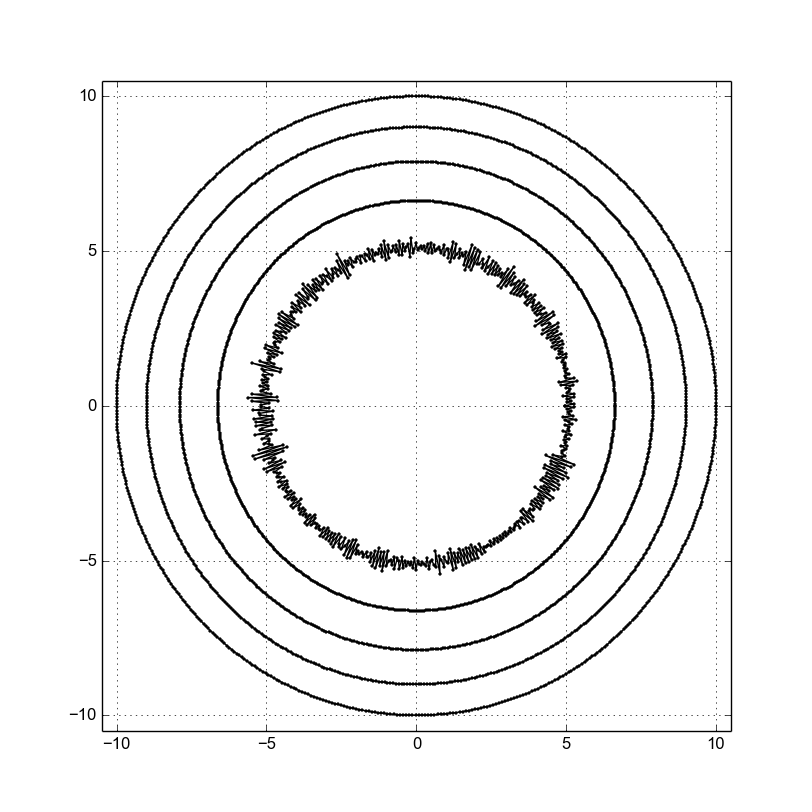}
   \includegraphics[scale=0.3]{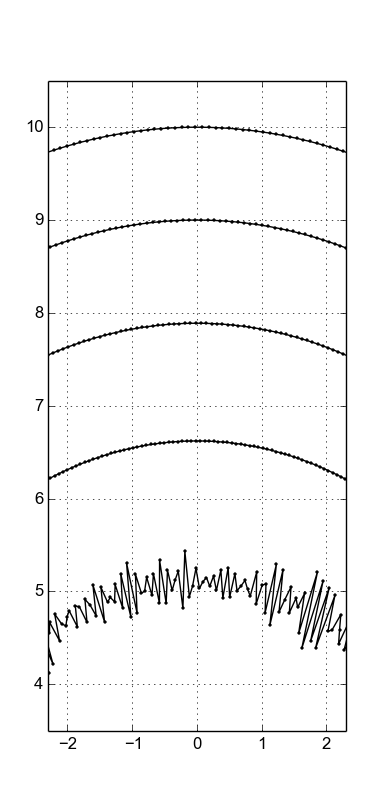}
  \caption{Numerical instability in a straightforward discretization of MCF for a circle (left), and detail (right). \label{fig:instability}}
\end{figure}

\begin{figure}
  \centering
  \includegraphics[scale=0.38]{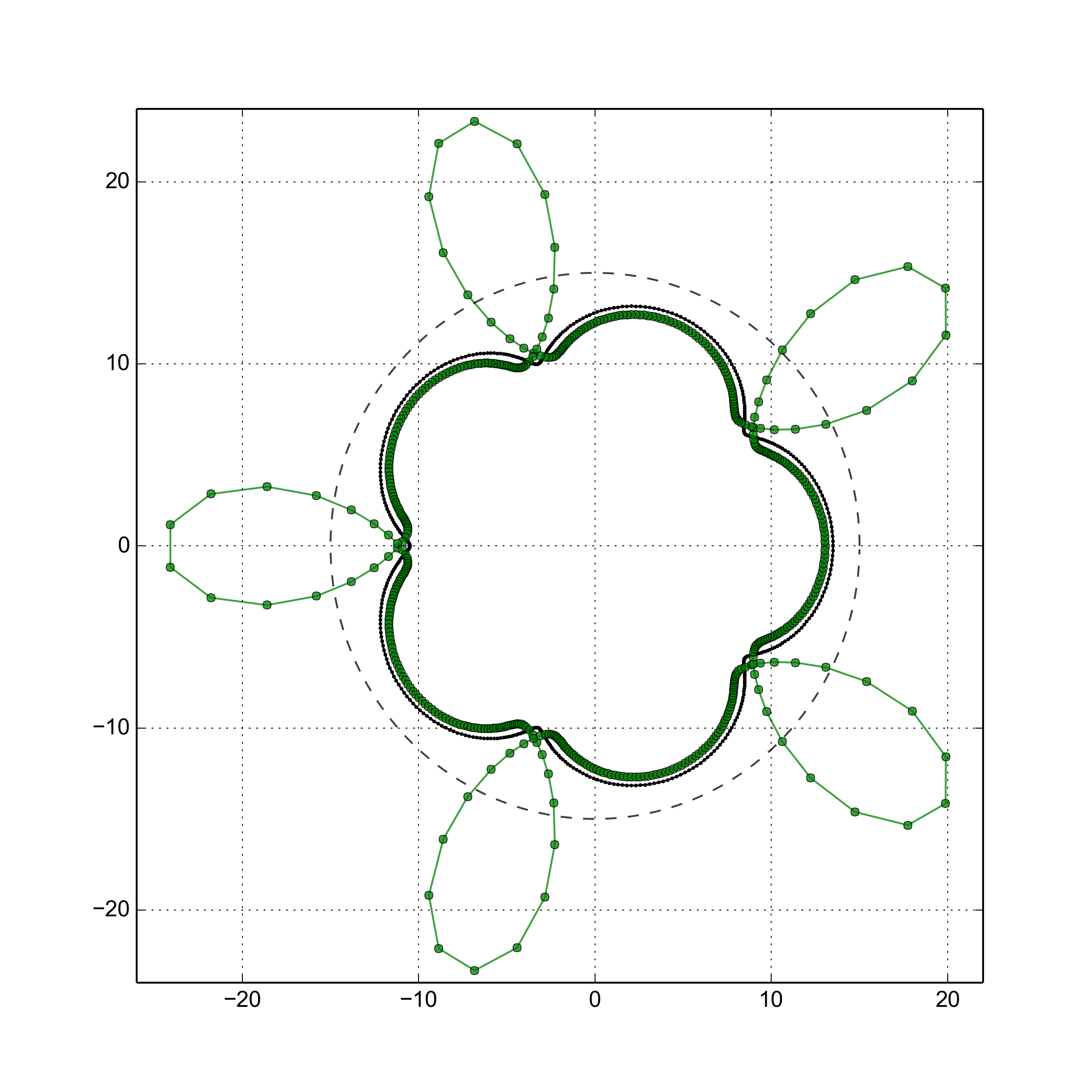}
   \includegraphics[scale=0.38]{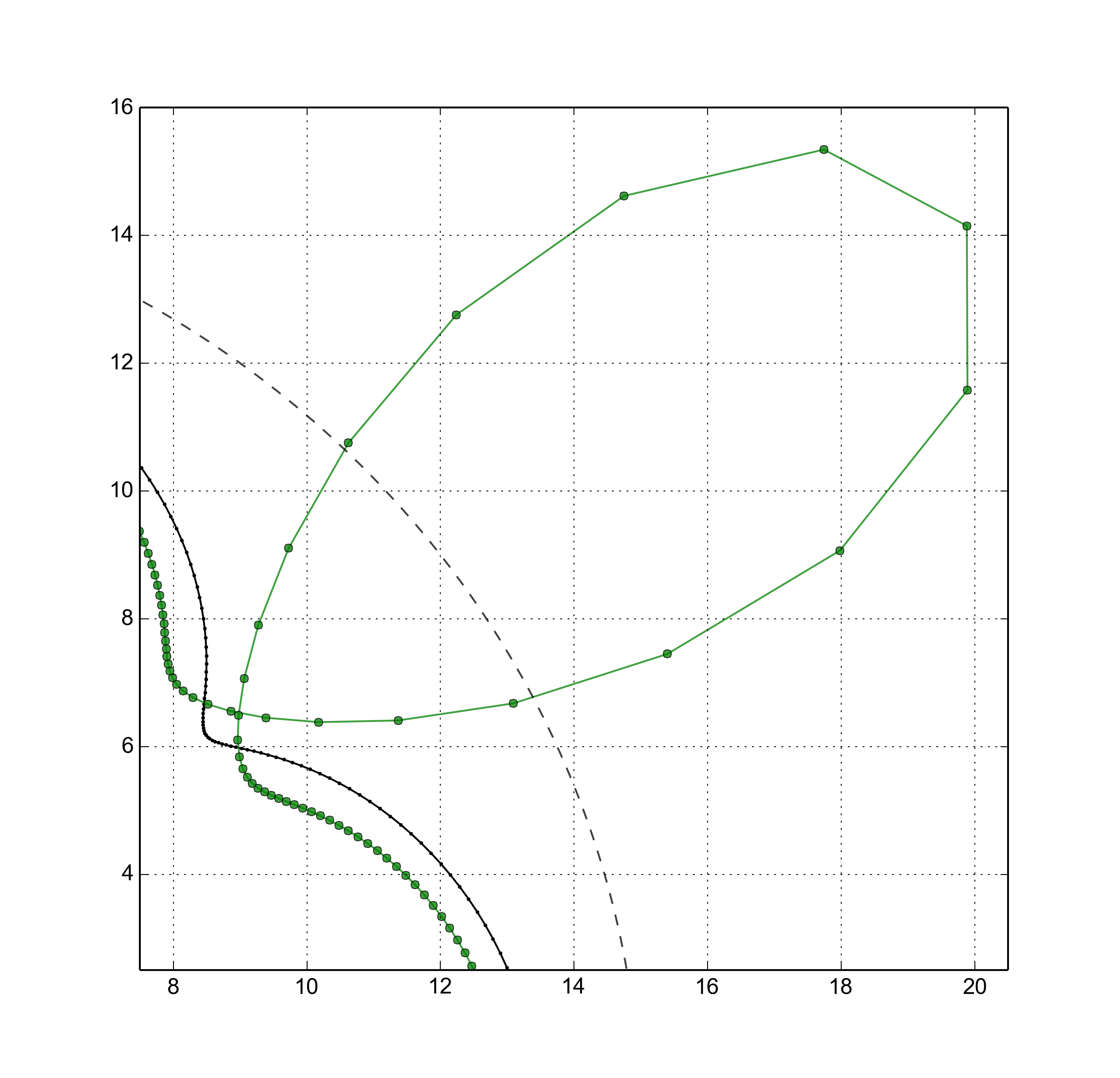}
  \caption{Loop formation in numerical discretization by standard Euler method of equation~(\ref{eq:MCF}) for a cycloid (left), and detail (right). \label{fig:loopform}}
\end{figure}

The disadvantage of Lagrangian methods are discussed in many works. Two problems arise commonly when straightforward discretization of equation (~\ref{eq:MCF}) is performed: (i) a numerical instability as shown in Figure~\ref{fig:instability} and (ii) loop
formation which contradicts the comparison principle (see~\cite{ecker,mantegazza}), as shown in Figure~\ref{fig:loopform} for a cycloid (solid line) and its first iteration result (polygonal curve) compared with a circumference (dashed curve).

Here we present a Lagrangian method for contour parametrization. To accomplish this application, we assume that a planar closed and differentiable curve is drawn on a 2D digital image. Then, we evolve the curve by mean curvature flow, but constraining the motion of the curve by the objects in the image. If only one object is initially inside the curve, as the curve shrinks, it will match parts of the boundary of the object. The main problem using MCF to recognize images is to couple the restriction and the flow, because it implies that the curve is not evolving uniformly. That is, not all points in the curve will move, even when the curvature at those points is different from zero. This kind of flow is called \emph{anisotropic}. To avoid certain numerical difficulties, our scheme considers a curve motion along tangential and normal directions, as in~\cite{kimura}. The results of~\cite{kimura} correspond to an unconstrained flow, so the stability and convergence results found there can not be compared with ours. Since the numerical procedure is based on certain MCF properties, we provide the necessary technical details about existence and uniqueness of solutions for MCF in Theorems~\ref{teo:linear} and~\ref{teo:exist}.

\begin{figure}
  \centering
  \includegraphics[scale=0.30]{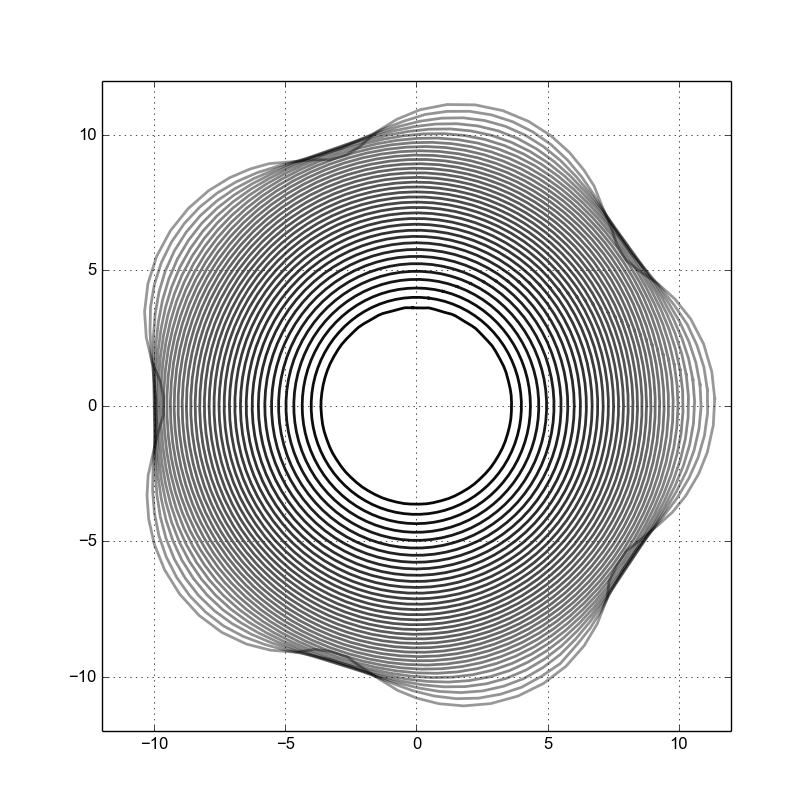}
   \includegraphics[scale=0.35]{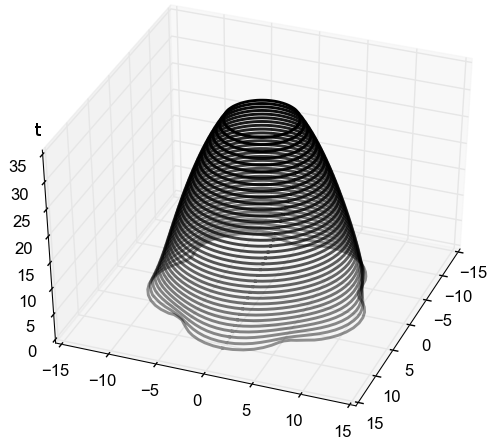}\\
   (a)\hspace{0.33\textwidth} (b)
  \caption{Numerical MCF of cycloid. (a) An evolving cycloid (outer curve) evolves by MCF, as time increase, it converges to a circumference. (b) 3D projection of (a), the time parameter $t$ is represented by the vertical axis. \label{fig:Num-examples1}}
\end{figure}

\begin{figure}
  \centering
  \includegraphics[scale=0.30]{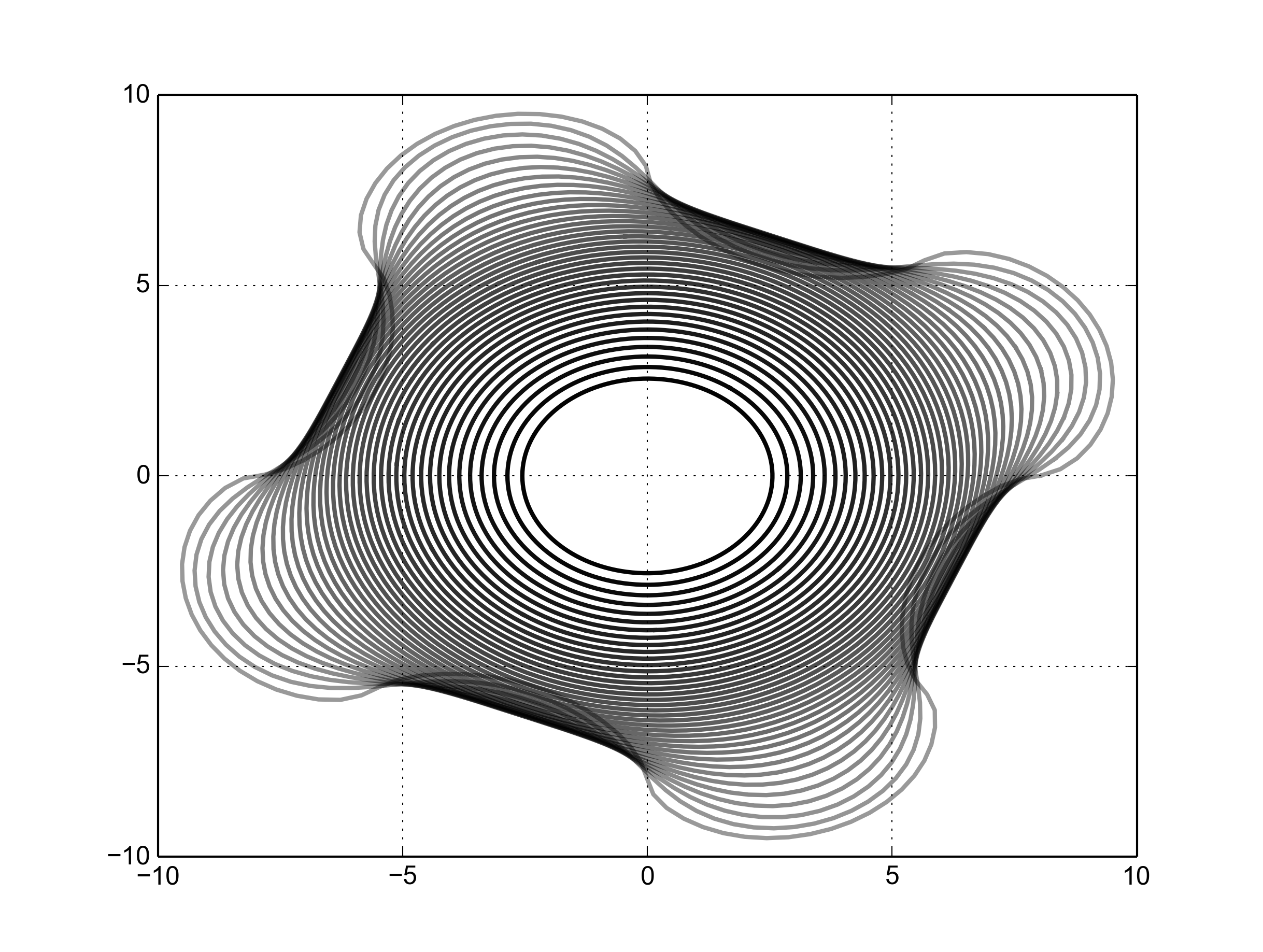}
   \includegraphics[scale=0.35]{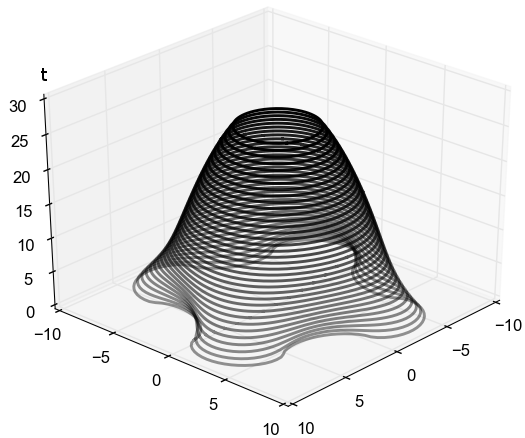}\\
   (a)\hspace{0.33\textwidth} (b)
  \caption{Numerical MCF of non-convex curve (I). (a) The outer curve is the initial non-convex curve, it deforms to a circumference and shrinks as time elapses. (b) 3D projection of (a), the time parameter $t$ is represented by the vertical axis.\label{fig:Num-examples2}}
\end{figure}

\begin{figure}
  \centering
  \includegraphics[scale=0.30]{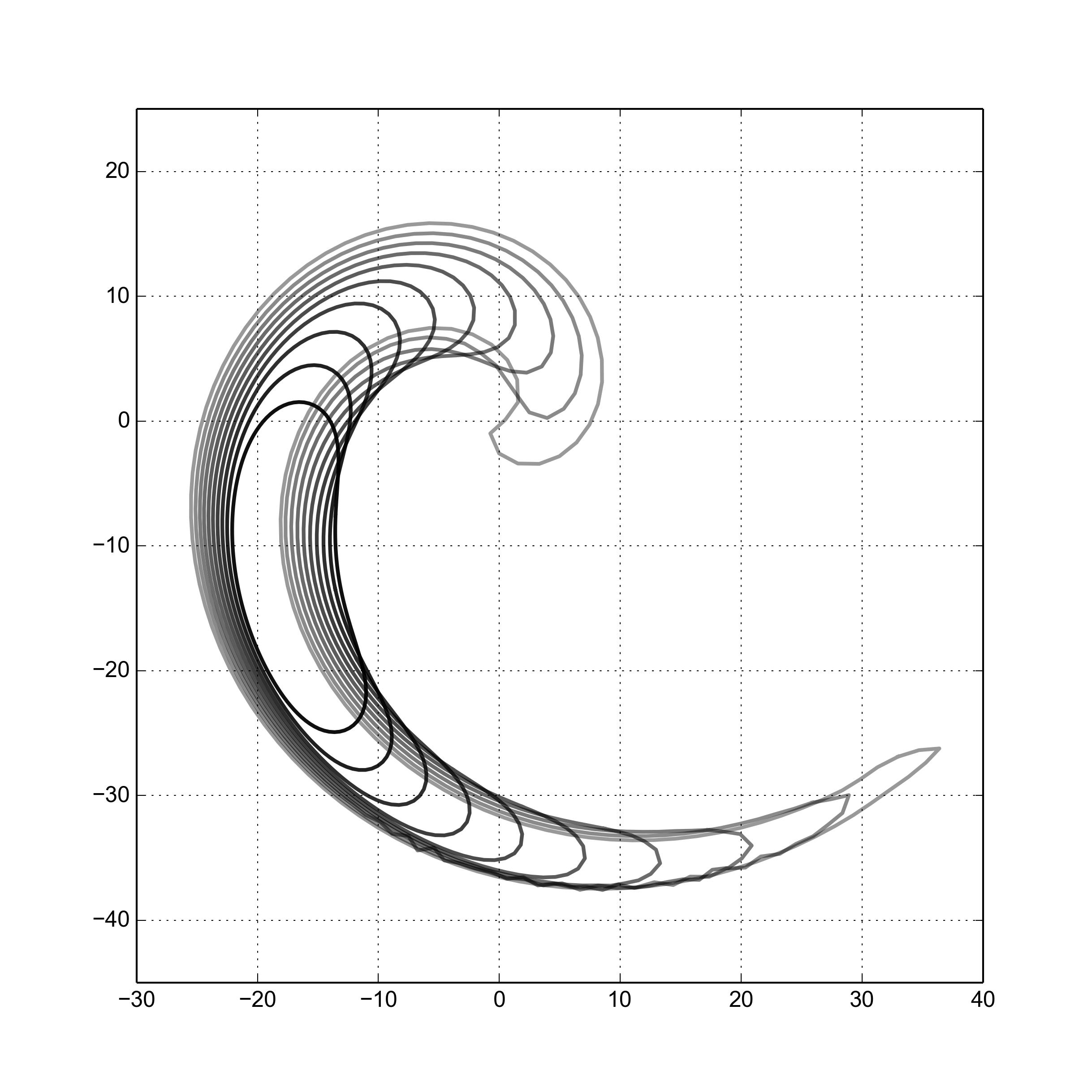}
   \includegraphics[scale=0.35]{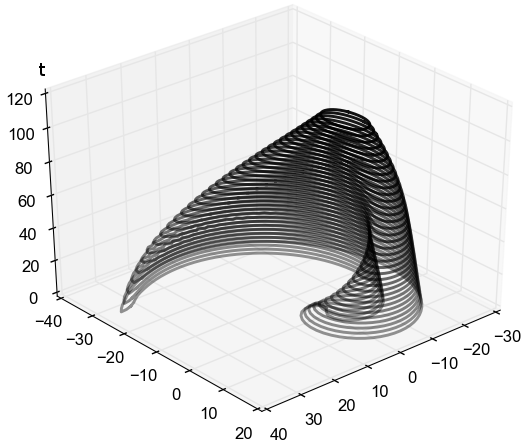}\\
   (a)\hspace{0.33\textwidth} (b)
  \caption{Numerical MCF of non-convex curve (II). (a) The outer curve is the initial non-convex curve, it shrinks as time elapses. (b) 3D projection of (a), the time parameter $t$ is represented by the vertical axis. \label{fig:Num-examples3}}
\end{figure}

Kimura developed a Lagrangian method which numerically reproduces the mean curvature flow for curves, based on a redistribution of points by reparametrization by arc length ~\cite{kimura}. That method represents an initial simply closed smooth curve $\varphi(s)$ by an ordered set of $N$ points. The position of  each point is updated to represent the curve after a time $\Delta t $. This method assumes a $\delta-$tubular neighborhood in which, for a small time interval $\Delta t$, the whole evolution is inside, and therefore may be parameterized as
\begin{gather*}
X(s,r;t)=\varphi(s)+r\mathbf{n}(s,t).
\end{gather*}
In this tubular neighborhood the map $X$, given by
\begin{gather*}
(s,r)\mapsto \{ x\in \mathbb{R}^{2} \, |\, {\rm dist}(x,\varphi(t))<\delta\}
\end{gather*}
is invertible. A correction term is then added to Euler's formula (see~\cite{kimura}). Figures~(\ref{fig:Num-examples1}-~\ref{fig:Num-examples3}) show numerical examples with Kimura's method, the evolution in time is represented by the z axis.\\

Our main theoretical result is the postulation of an {\it anisotropic} MCF, for which short time existence of solutions is shown in Theorem (\ref{teo:P1P2}). Furthermore, we detail a numerical scheme to approximate solutions and apply it to the task of contour parametrization.

\begin{figure}
  \centering
  \includegraphics[scale=0.28]{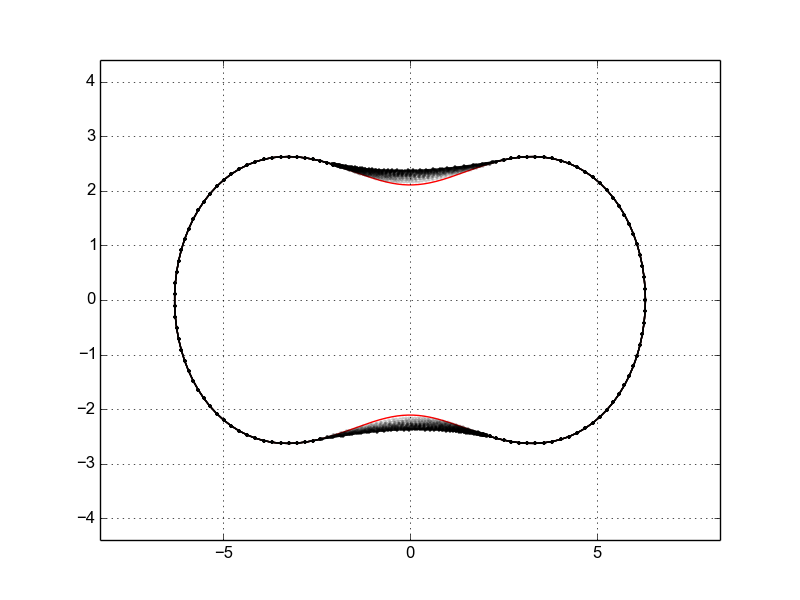}
   \includegraphics[scale=0.28]{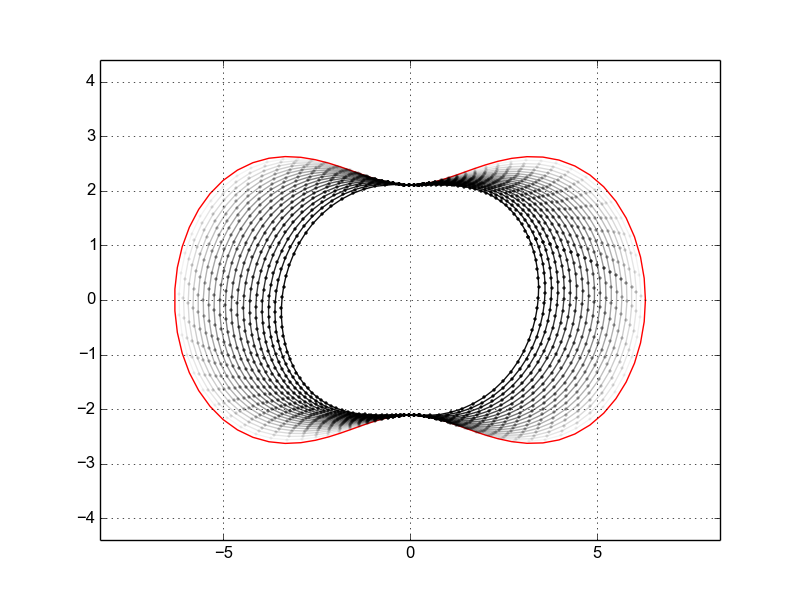}
   \includegraphics[scale=0.28]{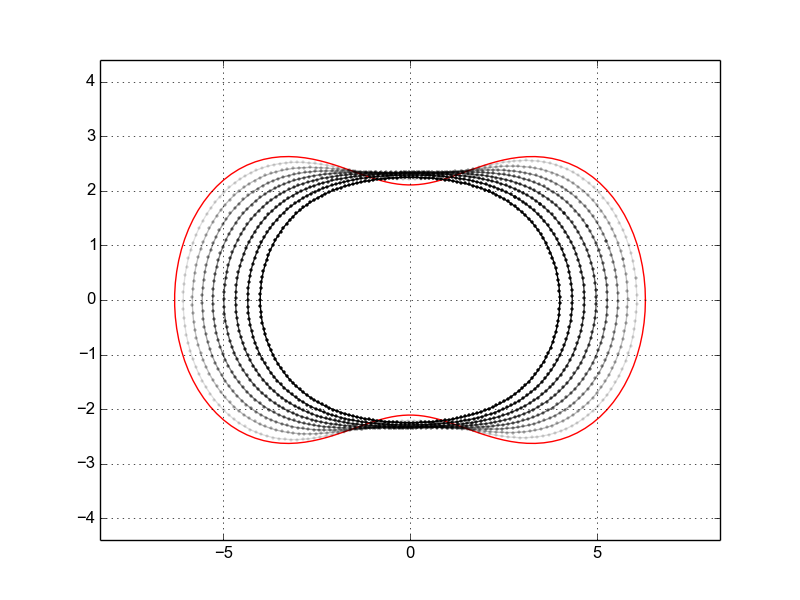}
  \caption{A non convex red curve evolves by anisotropic MCF (\emph{left and center}) and MCF (\emph{right}). In the figure on the left the evolution is obtained by replacing in the MCF equation the curvature $\kappa$ of the curve with $\min (\kappa,0)$. In the figure in the center, we replace $\kappa$ by $\max(\kappa,0)$. \label{fig:anisotropic}}
\end{figure}

Many varieties of an anisotropic MCF can be proposed, Figure (\ref{fig:anisotropic}) compares the MCF with two versions of anisotropic MCF. Our evolution scheme is an adaptation for contour recognition of Kimura's planar curve evolution. We approach the stability of our scheme through von Neumann's analysis. In our analysis, the parameters in error propagators (see Equation (\ref{ec:ampl_E})) are time dependent. Then, the stability condition cannot be determined for all time but only for the next time step. We establish our main stability criterion in Proposition~\ref{prop:criteria1}.\\

\begin{figure}
\centering
\includegraphics[scale=0.5]{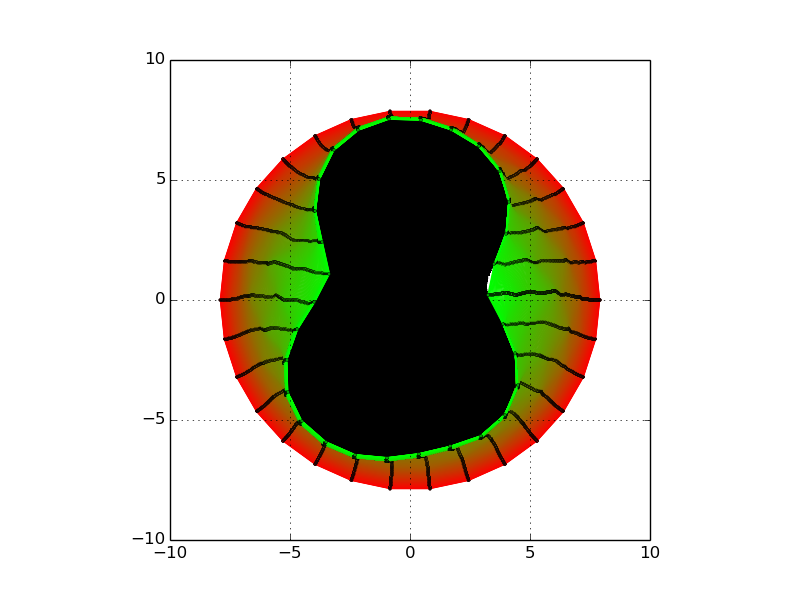}
\caption{A non-convex contour parametrization example of our technique. An initial red circumference is evolving by our AMCF, the enclosed black region constrains its motion (see section (\ref{sec:Implementation})). The color gradient from red to green represents the initial and final curve position. A high contrast between red and green represents a fast matching.\label{fig:nonconvexmatching}}
\end{figure}

In contrast to previous approaches, our method is a Lagrangian scheme whose main features are: (I) Estimation of curvature bounds are not required. (II) It is formulated as a Poisson problem with a boundary condition given implicitly. This condition links Poisson's problem with the MCF. Our proposed Poisson problem couples MCF and the field due to a point charge distribution. (III) Finding a solution requires to solve a $2N$ linear system. (IV) When it is implemented to a contour parametrization problem, only the pixel value at each point constraints the motion. (V) These constraints are handled through the source function and the boundary condition in Poisson's equation. Our numerical experiments reveal that this Poisson formulation can match some non-convex shapes perfectly. In general, matching non-convex shapes successfully depends on the charge distribution. Our scheme takes the numerical approximation for curvature and an equidistant distribution of points from~\cite{kimura}, and through our proposed Poisson formulation we can handle these constraints. In Figure~\ref{fig:nonconvexmatching}, we present a contour matching with an initial circumference as the evolving curve. The color gradient advances with time from red to green. A high contrast between red and green represents a fast matching, and a smooth color gradient indicates regions which require more iterations to reach the shape boundary in black. 

\section*{Outline}
This work is presented as follows: \S 1  the proof of existence and uniqueness of solutions for MCF and some of its properties, included for completeness sake; \S 2  the presentation of our conditioned or anisotropic MCF and formulas for its solutions; \S 3 the numerical implementation of our anisotropic flow; \S 4 the implementation of the anisotropic MCF for contour parametrization; \S 5 details the stability properties of our numerical scheme, and \S 6 contains our conclusions and suggestions for future directions.

\section{Short-time Existence and Uniqueness of Solutions}
The main result of this section proves existence and uniqueness of solutions to an isotropic version of MCF. The proof includes new technical details, in particular bounds on the displacement of the immersed curve (equation~(\ref{ec:reparamTubularN}) in Lemma~\ref{lemma:forCurves}), which are essential for the computations in the associated numerical problem in section~\ref{sec:numerical}.\\

This proof requires to rewrite equation~(\ref{eq:MCF}) in terms of a reparametrization of the evolving manifold. Recall the following two theorems.

\begin{theorem}\label{th:rewriteeq}(Xi-Ping Zhu~\cite{zhulectures})
The mean curvature flow equation for a hypersurface $M\subset \mathbb{R}$ with metric tensor $g$ is equivalent to
\begin{gather}
\frac{\partial X}{\partial t}=\Delta_{g} X
\end{gather}
\end{theorem}

For the reader's convenience we include a proof here.
\begin{proof}
Let $A$ be the second fundamental form of a manifold $M$ and $\nabla $ the induced connection by the ambient manifold connection $\bar{\nabla}$. In local coordinates the right side of mean curvature flow equation expands to
\begin{eqnarray}
\displaystyle H\,\mathbf{\hat{n}}\,&=&\sum _{ij}g^{ij}A_{ij}\mathbf{\hat{n}}=\sum _{ij}g^{ij}\mbox{\scalebox{2}{$\langle$}}\bar{\nabla}_{\frac{\partial X}{\partial x^{i}}}\, \frac{\partial X}{\partial x{j}},\mathbf{\hat{n}} \mbox{\scalebox{2}{$\rangle$}}\, \mathbf{\hat{n}}=\sum _{ij}g^{ij}\left( \bar{\nabla}_{\frac{\partial X}{\partial x^{i}}}\, \frac{\partial X}{\partial x{j}}-\left( \bar{\nabla}_{\frac{\partial X}{\partial x^{j}}} \frac{\partial X}{\partial x^{i}}\right)^{\top} \right) \nonumber
\end{eqnarray}
We focus on the individual components and develop further:
\begin{eqnarray}
\displaystyle H\,\mathbf{\hat{n}}\, ^{\alpha}\,&=& \sum_{i,j}g^{ij}\left( \frac{\partial ^{2}X}{\partial x^{i}\partial x^{j}} ^{\alpha}\,-\, \sum _{k,l,\beta}g^{kl} \frac{\partial X}{\partial x^{k}}^{\alpha}\,\frac{\partial X}{\partial x^{l}}^{\beta}\,\frac{\partial ^{2}X}{\partial x^{i}\partial x^{j}} ^{\beta}\right)\nonumber\\
&=&\sum_{i,j}g^{ij}\frac{\partial ^{2}X}{\partial x^{i}\partial x^{j}}^{\alpha}\,-\,\sum_{i,j}g^{ij}\sum_{k,l,\beta}g^{kl}\, \frac{\partial X}{\partial x^{k}}^{\alpha}\,\frac{\partial X}{\partial x^{l}}^{\beta}\,\frac{\partial ^{2}X_{l}}{\partial x^{i}\partial x^{j}}^{\beta}\nonumber\\
&=&\sum_{i,j}g^{ij}\frac{\partial ^{2}X}{\partial x^{i}\partial x^{j}}^{\alpha}\,-\,\sum_{i,j}g^{ij}\sum_{k,l,\beta}\frac{\partial X}{\partial x^{k}}^{\alpha}\,g^{kl}\left( \frac{\partial X}{\partial x^{l}}^{\beta}\,\frac{\partial ^{2}X}{\partial x^{i}\partial x^{j}}^{\beta}\,+\,\frac{1}{2}\,\frac{\partial ^{2}X}{\partial x^{i}\partial x^{l}}^{\beta}\frac{\partial X}{\partial x^{j}}^{\beta}\right. \nonumber \\
& &\left. \,-\,\frac{1}{2}\,\frac{\partial ^{2}X}{\partial x^{i}\partial x^{l}}^{\beta}\,\frac{\partial X}{\partial x^{j}}^{\beta}\,+\,\frac{1}{2}\,\frac{\partial ^{2}X}{\partial x^{j}\partial x^{l}}^{\beta}\,\frac{\partial X}{\partial x^{i}}^{\beta}\,\,-\,\frac{1}{2}\,\frac{\partial ^{2}X}{\partial x^{j}\partial x^{l}}^{\beta}\,\frac{\partial X}{\partial x^{i}}^{\beta}\right)\nonumber\\
&=&\sum_{i,j}g^{ij}\frac{\partial ^{2}X}{\partial x^{i}\partial x^{j}}^{\alpha}\,-\,\sum_{i,j}g^{ij}\sum_{k,l}\frac{\partial X}{\partial x^{k}}^{\alpha}\,g^{kl}\frac{1}{2}\left( \frac{\partial}{\partial x^{i}}\mbox{\scalebox{2}{$\langle$}} \frac{\partial X}{\partial x^{j}},\frac{\partial X}{\partial x^{l}}\mbox{\scalebox{2}{$\rangle$}}\,\right. \nonumber\\
& &\left. +\, \frac{\partial }{\partial x^{j}}\mbox{\scalebox{2}{$\langle$}}\frac{\partial X}{\partial x^{i}},\frac{\partial X}{\partial x^{l}} \mbox{\scalebox{2}{$\rangle$}}\,-\,\frac{\partial}{\partial x^{l}}\mbox{\scalebox{2}{$\langle$}}\frac{\partial X}{\partial x^{i}},\frac{\partial X}{\partial x^{j}} \mbox{\scalebox{2}{$\rangle$}}\right)\nonumber\\
&=&\sum_{i,j}g^{ij}\frac{\partial ^{2}X}{\partial x^{i}\partial x^{j}}^{\alpha}\,-\,\sum_{i,j}g^{ij}\sum_{k}\frac{\partial X}{\partial x^{k}}^{\alpha}\,\frac{1}{2}\sum_{l}g^{kl}\left( \frac{\partial}{\partial x^{i}}g_{jl}\,+\, \frac{\partial }{\partial x^{j}}g_{il}\,-\,\frac{\partial}{\partial x^{l}}g_{ij}\right)\nonumber\\
&=&\sum_{i,j}g^{ij}\frac{\partial ^{2}X}{\partial x^{i}\partial x^{j}}^{\alpha}\,-\,\sum_{i,j}g^{ij}\sum_{k}\frac{\partial X}{\partial x^{k}}^{\alpha}\,\Gamma_{ij}^{k}\nonumber\\
&=&\sum_{i,j}g^{ij}\left( \frac{\partial ^{2}X}{\partial x^{i}\partial x^{j}}^{\alpha}\,-\ \sum_{k}\frac{\partial X}{\partial x^{k}}^{\alpha}\,\Gamma_{ij}^{k}\right)=\sum_{i,j}g^{ij}\left( \nabla_{\frac{\partial X}{\partial x^{i}}}\nabla_{\frac{\partial X}{\partial x^{j}}}X\right) ^{\alpha}=\Delta_{g} X^{\alpha}\nonumber
\end{eqnarray}\end{proof}

\begin{theorem}\label{teo:invariance} (Mantegazza~\cite{mantegazza})
Let $M^{n}$ be a compact submanifold of $\mathbb{R}^{n+1}$ with induced metric $g$ and $\tilde{X}:M\times [0,t)\rightarrow \mathbb{R}^{n+1}$ be an immersion which satisfies at every point $\mathbf{p}$ in $M$ and every time $t$ in $[0,\tau)$
\begin{equation}\label{ec:invariance}
\begin{split}
\frac{\partial \tilde{X}}{\partial t}(\mathbf{p},t)&=\Delta_{g}\tilde{X}(\mathbf{p},t)+F(\mathbf{p},t)\\
\tilde{X}(\cdot,0)&=\tilde{X}_{0}(\cdot)
\end{split}
\end{equation}
here $F(\mathbf{p},t)$ is in $d\tilde{X}_{t}|_{(\mathbf{p},t)}(T_{\mathbf{p}}M)$. Then, there exist a family of diffeomorphisms $\varphi$ such that $\tilde{X}\circ\varphi$ satisfies equation~(\ref{eq:MCF}). Conversely, given an immersion $X:M\times [0,\tau)\rightarrow\mathbb{R}^{n+1}$ and a reparametrization $\varphi$ such that $X\circ\varphi$ is a mean curvature flow, then there is a field $F$ in $dX_{t}|_{(\mathbf{p},t)}(T_{\mathbf{p}}M)$ such that satisfies equation~(\ref{ec:invariance}).
\end{theorem}

In the following, we restrict the study to plane curves, compare with~\cite{GerhardPolden1999} and~\cite{mantegazza}. Set a regular smooth closed plane curve $\gamma_{0}$ with unitary tangent vector $T_{0}$, normal vector $\mathbf{\hat{n}}_{0}$ and curvature $\kappa_{0}$.\\

The proof of short-time existence and uniqueness of solutions for the MCF is divided in three steps. We include these details here because they will be relevant to our main result in the analysis of solutions to {\it anisotropic} MCF in Theorem (\ref{teo:P1P2}) below.

\begin{quote}
\begin{itemize}
\item[{\bfseries Step 1}] Reparametrize the deformations and rewrite the evolution in terms of a scalar quasilinear parabolic problem.
\item[{\bfseries Step 2}] Prove existence and uniqueness of solutions for the linearized problem.
\item[{\bfseries Step 3}] Extend Step 2 to the quasilinear case using the inverse function theorem for Banach Spaces~\cite[Sec. 4.13]{zeidlerFPT}. 
\end{itemize}\end{quote}

\subsection*{Step 1}
First, we reparametrize the deformation of the curve in the normal direction for a small time interval.
\begin{lemma}\label{lemma:forCurves}
Let $\epsilon >0$ be small enough such that for $t$ in $[0,\epsilon)$ the evolution by mean curvature of $\gamma_{0}$ is within a tubular neighborhood. If the deformations are parameterized by
\begin{gather}\label{ec:reparamTubularN}
\gamma (s,t)=\gamma_{0}(s)+f(s,t)\mathbf{\hat{n}}_{0}
\end{gather}
Then, the evolution equation for~$(\mathrm{\ref{ec:reparamTubularN}})$ is given by
\begin{gather}\label{ec:dfdt}
\frac{\partial f}{\partial t}=\frac{1-\kappa_{0}}{\sqrt{(1-\kappa_{0}f)^{2}+(f^{\prime})^{2}}}\kappa(t)
\end{gather}
\end{lemma}
\begin{proof}
Using equation~(\ref{ec:reparamTubularN}) we compute the unitary tangent vector $T$, the normal vector $\mathbf{\hat{n}}$ and the curvature $\kappa$ of $\gamma(s,t)$
\begin{equation}
\begin{split}
T&=\frac{\gamma ^{\prime}}{||\gamma ^{\prime}||}=\frac{\gamma ^{\prime}_{0}+f^{\prime}\mathbf{\hat{n}}_{0}-\kappa_{0} fT_{0}}{||\gamma ^{\prime}_{0}+f^{\prime}\mathbf{\hat{n}}_{0}-\kappa_{0} fT_{0}||}=\frac{(1-\kappa_{0} f)T_{0}+ f^{\prime}\mathbf{\hat{n}}_{0}}{||(1-\kappa _{0}f)T_{0}+ f^{\prime}\mathbf{\hat{n}}_{0}||}\\
&=\frac{(1-\kappa_{0} f)T_{0}+f^{\prime}\mathbf{\hat{n}}_{0}}{\sqrt{ (1-\kappa_{0} f)^{2}+(f^{\prime})^{2}}}
\end{split}
\end{equation}

\begin{equation}
\begin{split}
\mathbf{\hat{n}}&=\frac{\mathbf{\hat{n}}_{0}-\frac{f^{\prime}}{||\gamma^{\prime}||}\,\mathbf{T}}{||\mathbf{\hat{n}}_{0}-\frac{f^{\prime}}{||\gamma^{\prime}||}\,\mathbf{T}||}=\frac{\left(1-\frac{(f^{\prime})^{2}}{||\gamma ^{\prime}||^{2}}\right)\mathbf{\hat{n}}_{0}-\frac{f^{\prime}(1-\kappa_{0}f)}{||\gamma^{\prime}||^{2}}\mathbf{T}_{0}}{||\mathbf{\hat{n}}_{0}-f^{\prime}\,\mathbf{T}||}\\
&=\frac{\left(1-\frac{(f^{\prime})^{2}}{||\gamma ^{\prime}||^{2}}\right)\mathbf{\hat{n}}_{0}-\frac{f^{\prime}(1-\kappa_{0}f)}{||\gamma^{\prime}||^{2}}\mathbf{T}_{0}}{\sqrt{(1-\frac{(f^{\prime})^{2}}{||\gamma^{\prime}||^{2}})^{2}+(f^{\prime})^{2}\frac{(1-\kappa_{0}f)^{2}}{||\gamma^{\prime}||^{4}}}}\\
&=\frac{(1-\kappa_{0}f)^{2}\mathbf{\hat{n}}_{0}-f^{\prime}(1-\kappa _{0}f)\mathbf{T}_{0}}{(1-\kappa _{0}f)\sqrt{(1-\kappa_{0}f)+(f^{\prime})^{2}}}=\frac{(1-\kappa_{0}f)\mathbf{\hat{n}}_{0}+f^{\prime}\,\mathbf{T}_{0}}{\sqrt{(1-\kappa_{0}f)^{2}+(f^{\prime})^{2}}}
\end{split}
\end{equation}

\begin{equation}\label{ec:kurvt}
\begin{split}
\kappa (t)&=\frac{|| \mathbf{T}^{\prime} \wedge \mathbf{T}^{\prime\prime} ||}{||\gamma ^{\prime}||^{3}}\\
&=\frac{||\{ (1-\kappa_{0}f)\mathbf{T}_{0}+f^{\prime}\mathbf{\hat{n}}_{0}\}\wedge \{ (-\kappa ^{\prime}_{0}-2\kappa_{0}f^{\prime})\mathbf{T}_{0}+(\kappa_{0}(1-\kappa_{0}f)+f^{\prime\prime})\mathbf{\hat{n}}_{0}\}||}{||\gamma ^{\prime}||^{3}}
\\
&=\frac{(1-\kappa_{0}f)f^{\prime\prime}+2\kappa _{0}(f^{\prime})^{2}+\kappa_{0}^{\prime}ff^{\prime}-2\kappa_{0}^{2}f+\kappa_{0}^{3}f^{2}}{\left( (1-\kappa_{0}f)^{2}+(f^{\prime})^{2}\right)^{3/2}}
\end{split}
\end{equation}
As $\gamma(s,t)$ solves~(\ref{eq:MCF}), we differentiate equation~(\ref{ec:reparamTubularN}) and take the component along $\mathbf{\hat{n}}_{0}$:
\begin{gather*}
\frac{\partial f}{\partial t} \langle \mathbf{\hat{n}}_{0},\mathbf{\hat{n}}_{0}\rangle=\kappa(t)\,\,\langle \mathbf{\hat{n}},\mathbf{\hat{n}}_{0}\rangle\mbox{ .}
\end{gather*}
Substituting and solving for $\partial f/\partial t$ we obtain
\begin{gather*}
\frac{\partial f}{\partial t}=\frac{1-\kappa_{0}}{\sqrt{(1-\kappa_{0}f)^{2}+(f^{\prime})^{2}}}\,\kappa(t)\qedhere
\end{gather*}\end{proof}

\subsection*{Step 2}
Denote by $W_{2}^{k}$ the set of $L_{2}$ functions over $\gamma_{0}$ whose $k$ derivatives are also in $L_{2}$. Let $u$ an $v$ be functions in $W_{2}^{2}$, consider the inner product
\begin{gather*}
\langle \, u\,,\,v\,\rangle_{W_{2}^{k}}=\int_{\gamma_{0}}\,\,\sum_{i,j\leq k}\, \frac{\partial^{i} u}{\partial x^{i}}\frac{\partial^{j} v}{ \,\partial x^{j}} \,dx,
\end{gather*}
and denote by $W_{2}^{\alpha,\beta}$ the set of $L_{2}$ functions over $\gamma_{0}\times(0,\infty)$ whose $\alpha$ spatial derivatives and $\beta$ time derivatives are also in $L_{2}$.
We now state the linearized problem:

\begin{theorem}\label{teo:linear} 
Let $u_{0}$ be a $W_{2}^{\alpha +1}$ function over $\gamma_{0}$ , $h\,\in\,W_{2}^{\alpha,\beta}$, and $\mathcal{L}_{0}$ be a linear parabolic operator. The problem
\begin{equation*}
\mbox{({\bfseries LP})}\;\begin{cases}
&\mathcal{L}_{0}(u)=h\\
&u(\cdot,0)= u_{0}
\end{cases}
\end{equation*}
has a solution.
\end{theorem}

Before proceeding to the proof of theorem~\ref{teo:linear}, we will need the following:

\begin{lemma}\label{lema:garding} (Garding's inequality~\cite{zeidler89}) Let $k>0$ an integer, suppose that the application of a linear differential operator $\mathcal{A}$ over a $2k$-times differentiable function $f$ over $M^{n}\times[0,T)$ is given by
\begin{gather*}
A(f(x)):=\sum_{\alpha,\beta=1}^{n}\partial_{\alpha}( a_{\alpha\beta}\partial_{\beta}f)+\sum_{\alpha=1}^{n}\partial_{\alpha}(b_{\alpha} f)+cf.
\end{gather*}
Let $\xi$ be in $\mathbf{R}^{k}$ and $\lambda$, $\mu$ be positive constants. Suppose that the coefficients $a_{\alpha\beta}$ are bounded and satisfy
\begin{gather}\label{ec:condalpha}
\lambda ||\xi ||^{2}\leq \sum_{\alpha,\beta=1}^{n}a_{\alpha \beta}\,\xi^{\alpha}\xi^{\beta}.
\end{gather}
Then, if $u(x)$ is in $W_{2}^{k}(M)$, there are constants $c_{1}>0$ and $c_{2}$ such that
\begin{gather*}
\int_{M}\left( \sum_{\alpha,\beta=1}^{n}a_{\alpha \beta}\partial_{\alpha}u\partial_{\beta} u\,-\sum_{\alpha=1}^{n}b_{\alpha}u\partial_{\alpha}u+cu^{2} \right)\,dx\geq c_{1}||u||^{2}_{W^{k}_{2}(M)}-c_{2}||u||^{2}_{L^{2}(M)}
\end{gather*}
\end{lemma}

\begin{lemma}\label{lema:laxM}(Modified Lax-Milgram theorem~\cite{wu2006elliptic}) Let $H$ be a Hilbert space, $V\subset
H$ a dense subspace and $a(u,v)$ a bilinear form in $H\times V$ such that
\begin{itemize}
\item[(1)] Let $u\in H$, $v\in V$ and some $C>0$
\begin{gather*}
|a(u,v)|\leq C\,||u||_{H}\,||v||_{V}.
\end{gather*}
\item[(2)] For every $v\in V$ and some $\delta >0$\\
\begin{gather*}
|a(v,v)|\geq \delta \, ||v||^{2}_{H}
\end{gather*}
\end{itemize}
Then for every bounded linear operator $F$ in $H$ and $v$ in $V$, there exists $u$ in $H$ such that $F(v)=a(u,v)$.
\end{lemma}

Let $f,\, g$ in $C^{\infty}_{0}(\gamma_{0}\times[0,\infty))$, we recall the following definitions from~\cite{GerhardPolden1999}:
\begin{gather*}
\langle f,g\rangle_{LL_{m}}=\int_{0}^{\infty}e^{-mt}\,f\,g\, dt\\
\langle f,g \rangle_{LW_{m}^{k}}=\int_{0}^{\infty}e^{-mt}\langle f,g \rangle_{W_{2}^{k}}\, dt\,.
\end{gather*}

Here $LL$ and $LW$ are the Hilbert spaces resulting from the completion of the space of $C^{\infty}_{0}(\gamma_{0}\times[0,\infty))$ functions with compact support using the norms above, in that order.\\

Let $V$ be the space of $C^{\infty}$ functions over $(\gamma_{0}\times[0,\infty))$ such that $V(\cdot,t)=0$ for small and large values of $t$, and let $WW^{k}$ be the completion of $V$ with the norm associated to the following inner product
\begin{gather*}
\langle\, f\,,\,g\,\rangle_{WW^{k}}=\langle f\,,\,g \rangle_{LW_{m}^{k}}+\langle f_{t}\,,\,g_{t} \rangle_{LL_{m}}.
\end{gather*}
Let
\begin{gather*}
P^{s}=\left\{ f:\gamma_{0}\times[0,\infty)\rightarrow\mathbb{R}\, :\,\,|\partial_{t}^{i}f|_{LW_{m}^{2}}<\infty,\, \forall i\leq s\right\}
\end{gather*}
endowed with
\begin{gather*}
\langle f,g\rangle_{P^{s}}=\sum_{i\leq s}\langle \partial_{t}^{i}f\,,\,\partial_{t}^{i}g\rangle_{LW_{m}^{2(s-i)} }.
\end{gather*}

From the theory of parabolic linear partial differential equations theory we invoke ~\cite[p. 351-352]{evans2010partial}.

\begin{definition}Let $Q:=\gamma_0 \times (0,\epsilon)$ and $\varphi$ be a $W_{2}^{1,1}$ function over $\gamma_{0}\times [0,\infty)$ which vanishes at $t=0$. We say that $u(x,t)$ in $WW$ is a \emph{weak solution} of the linearized problem ({\bfseries LP}) if
\begin{gather}\label{ec:weakSol}
\int_{Q} (u_{t}\varphi+a(x,t)u_{x}\varphi _{x} +b(x,t)u_{x}\varphi+c(x,t)u\varphi)\, dxdt=\int_{Q} h(x,t)\varphi\, dxdt
\end{gather}
\end{definition}
The following lemma will be useful to construct a bounded bilinear operator, and then we will be able to apply Lemma~\ref{lema:laxM}.
\begin{lemma}(See~\cite[Sec. 7.1.2]{evans2010partial} and \cite[Sec. 3.5.1]{wu2006elliptic})\label{lema:weak}
Let $m$ be a positive constant. An $u(x,t)$ in $WW$ is a weak solution of ({\bfseries LP}) if and only if 
\begin{gather}
\int_{Q} (u_{t}\varphi_{t}+a(x,t)u_{x}\varphi_{tx} +b(x,t)u_{x}\varphi_{t}+c(x,t)u\varphi_{t})e^{-mt}\, dxdt=
\int_{Q} h(x,t)\varphi_{t}e^{-mt}\, dxdt\;.\;\;\label{ec:debint}
\end{gather}
\end{lemma}

We can now begin the proof of theorem~\ref{teo:linear}
\begin{proof}[Proof of Theorem~\ref{teo:linear}]
First we will prove the existence of a weak solution.\\

Let $A:WW^{2}\times V\rightarrow \mathbb{R}$ be a bilinear form given by
\begin{gather}\label{ec:defAbil}
A(u,v)=\int_{0}^{\infty}\int e^{-mt}\left( u_{t}v_{t}+a(x,t)u_{x}v_{tx}+b(x,t)u_{x}v_{t}+c(x,t)uv_{t} \right)\, dxdt
\end{gather}
here $m$ is a constant, $u\in WW^{2}$, $v\in V $, and $a,b,c$ are the coefficients of the linearized problem ({\bfseries LP}). Let $F(v)$ be a linear operator over $V $ given by
\begin{gather}\label{ec:defFv}
F(v):=\int_{0}^{\infty}\int h\,v\, e^{-mt}\, dxdt
\end{gather}
We will prove that $A$ satisfies the  hypotheses of lemma~\ref{lema:laxM}: 
\begin{itemize}
 \item Condition (1). The coefficients $a,\,b,\,c$ and $h$ are bounded because they were found by linearization of a bounded function in a tubular neighborhood, see Equation~\ref{ec:dfdt} in Lemma~\ref{lemma:forCurves}. Let $C$ be the maximum of the bounds for $a,\,b$ and $c$. Then,
\begin{eqnarray}
A(u,v)&\leq& \int_{0}^{\infty}\int e^{-mt}\left(u_{t}v_{t}+C(u_{x}v_{tx}+u_{x}v_{t}+uv_{t})\right)\,dxdt\nonumber\\
&\leq&\max(1,C)\int_{0}^{\infty}\int e^{-mt}\left(u_{t}v_{t}+u_{x}v_{tx}+u_{x}v_{t}+uv_{t} \right)\,dxdt\nonumber
\end{eqnarray}
Noticing that this double integral corresponds to  the inner product of $u$ and $v$ in $P^{1}_{2}$. Then,
\begin{eqnarray}
|A(u,v)|&\leq& \max(1,C)\,| \langle\, u\,,v\,\rangle_{P^{1}_{2}}|\leq \max(1,C)\,||u||_{P^{1}_{2}}\,||v||_{P^{1}_{2}}
\leq\max(1,C)||u||_{WW^{2}}\,|v|_{V}.\nonumber
\end{eqnarray}

\item To verify that condition (2) also holds observe that:
\begin{eqnarray}
A(v,v)&=&\int_{0}^{\infty}\int \left( v_{t}^{2}+au_{x}v_{tx}+b v_{x}v_{t}+cvv_{t}\right)e^{-mt}\, dx\, dt\nonumber\\
&=&\int_{0}^{\infty}\int v_{t}^{2}e^{-mt}\, dxdt\;+\;\int_{0}^{\infty}\int \left( au_{x}v_{xt}+b v_{x}v_{t}+cvv_{t} \right)e^{-mt}\, dx\, dt\nonumber\\
&=&||v_{t}||^{2}_{LL_{m}}\;+\;\int_{0}^{\infty}\int \left( av_{x}v_{xt}+b v_{x}v_{t}+cvv_{t} \right)e^{-mt}\, dx\, dt\nonumber
\end{eqnarray}
Integrate by parts to obtain:
\begin{eqnarray}
A(v,v)&=&||v_{t}||^{2}_{LL_{m}}\;+\;\frac{m}{2}\int_{0}^{\infty}\int e^{-mt}\left( av_{x}^{2}+bv_{x}v+cv^{2} \right)\,dxdt\nonumber\\
& &\;\;\;\;\;-\frac{1}{2}\int_{0}^{\infty}\int e^{-mt}\left( a_{t}v_{x}^{2}+b_{t}v_{x}v+c_{t}v^{2} \right)\,dxdt\,-\frac{1}{2}\int av_{x}^{2}|_{t=0}dx\nonumber
\end{eqnarray}
Let $\mathcal{S}$ be an upper bound for $a_{t},b_{t},c_{t}$. Using lemma~\ref{lema:garding} there exists constants $C_{1}$ and $C_{2}$ such that:\\
\begin{eqnarray}
A(v,v)&\geq &||v_{t}||^{2}_{LL_{m}}\;+\;\frac{m}{2}\int_{0}^{\infty}e^{-mt}C_{1}\, \,||v||_{W_{2}^{2}}^{2}\,dt-\;\frac{1}{2}\int_{0}^{\infty}e^{-mt}\mathcal{S}\,||v||_{W_{2}^{2}}^{2}\,dt\nonumber\\
& &-C_{2}\frac{m}{2}\int_{0}^{\infty}e^{-mt}||v||_{L_{2}}^{2}\,dt-\frac{C}{2}||v_{x}(0)||^{2}_{L_{2}}\nonumber\\
&=& ||v_{t}||^{2}_{LL_{m}}\;+\;\left( \frac{m\,C_{1}\,-\mathcal{S}}{2}\right) \,||v||_{LW}-\frac{C_{2}}{2}||v||_{L_{2}}^{2}-\frac{C}{2}||v_{x}(0)||^{2}_{L_{2}}\nonumber\\
&\geq &||v||_{WW}^{2}\nonumber
\end{eqnarray}
The last inequality follows from having properly chosen $m$.\\
\end{itemize}
From the above arguments we obtain that for all $v\in V$ there exists $u_{L}$ such that $F(v)=A(u_{L},v)$. Since  Equation~(\ref{ec:debint}) can be obtained from~(\ref{ec:defAbil}) and~(\ref{ec:defFv}), $u_{L} $ is therefore a weak solution.
\end{proof}

We proved the existence of weak solutions. However, the regularity of weak solutions for linear parabolic problems depends on the regularity of the data $h$ and $u_{0}$. Since the existence of solutions for linear parabolic problems is an auxiliary result in our work, we only refer to~\cite[Sec. 7.1.3]{evans2010partial} for the regularity extension.\\

The uniqueness of weak solutions for linear parabolic operators is also stated in~\cite{evans2010partial} by considering the difference  of two weak solutions $w:=u_{1}-u_{2}$. Since $w$ satisfies the problem $(\mathbf{LP})$ with $h=w(\cdot,0)=0$, from lemma~\ref{lema:garding} and substituting $w$ into~(\ref{ec:weakSol}, we obtain
\begin{gather}\label{eq:estimate}
\frac{d}{ds} |w(\cdot,s)|^{2}_{L_{2}}\leq C\,|w(\cdot,s)|_{L_{2}}^{2}
\end{gather}
Substituting the last inequality into the derivative of $|w|_{\,L_{2}}^{2}\,e^{-Cs}$, one obtains
\begin{gather*}
\frac{d}{ds}\left( |w(s)|_{L_{2}}^{2}\,e^{-Cs}\right)\leq 0
\end{gather*}
Since $|w(s)|_{L_{2}}^{2}\,e^{-Cs}$ is non negative, $|w(t)|_{L_{2}}^{2}\,e^{-Ct}-\cancel{|w(0)|_{L_{2}}^{2}}\leq 0$. Thus $w=0$ implies $\,u_{1}=u_{2}$. A similar procedure can be used to prove the uniqueness of solutions for quasilinear parabolic problems. The main hypothesis is an estimate similar to~(\ref{eq:estimate}), in particular for MCF see~\cite{Chen2006UniquenessAP}.

\subsection*{Step 3}
We now proceed to extend the linearized solution of normal deformation (Theorem~\ref{teo:linear}) to the non linear problem
\begin{theorem}\label{teo:exist} 
Given a quasilinear parabolic problem
\begin{equation*}
\begin{cases}
& \partial_{t}u -a(u,x,t)u_{xx}-b(u,x,t)u_{x}-c(u,x,t)u=h(u,x,t)\nonumber\\ 
& u(\cdot,0)= u_{0}.
\end{cases}
\end{equation*}
in $\gamma_{0}\times[0,\epsilon)$. Then, there exist $\epsilon>0$ such that this problem has a solution.
\end{theorem}
\begin{proof}
Define the map $\mathcal{L}:P^{2}\rightarrow W^{2,1}_{2}\times P^{1}$ by
\begin{gather*}
u \stackrel{\mathfrak{L}}{\mapsto} (u_{0},\mathcal{L}_{1}(u)).\nonumber
\end{gather*}
Notice that a function $u$ such that $\mathcal{L}(u)=0$ is a solution for the quasilinear problem. To this end, let $\mathcal{L}_{u_{0}}$ be the linearization of $\mathcal{L}_{1}$ around $u_{0}$, by Theorem~\ref{teo:linear} the problem
\begin{equation}\label{ec:aux}
\begin{cases}
& \mathcal{L}_{u_{0}}(w(x,t))= 0\\
& w(\cdot,0)= u_{0}
\end{cases}
\end{equation}
has a solution. Let $w^{t}:=w(\cdot,t)$ and $\{ w^{t}\}$ be a sequence converging to $u_{0}$ as $t$ tends to $0$ from above. Suppose that $w$ is a solution of~(\ref{ec:aux}), and let $a(u_{0};x,t),\,b(u_{0};x,t),\,c(u_{0};x,t)$ and $h(u_{0};x,t)$ be the coefficients of $\mathcal{L}_{u_{0}}$. Additionally, consider the linearization of $\mathcal{L}_{1}$ around $w$, let $a(w;x,t),\,b(w;x,t),\,c(w;x,t)$ and $h(w;x,t)$ be its linearized coefficients. From~(\ref{ec:aux}) we have
\begin{eqnarray}
\mathcal{L}_{1}(w)&=&\partial_{t}w-P(x,t,w,D_{x}w)\nonumber\\
&=&h(u_{0};x,t)+a(u_{0};x,t)w_{xx}+b(u_{0};x,t)w_{x}+c(u_{0};x,t)w-P(x,t,w,w_{x},w_{xx})\nonumber\\
&=& h(u_{0})u_{xx}-h(w)w_{xx} + a(u_{0})u_{x}-a(w)w_{x} +b(u_{0})u-b(w)w+c(u_{0})u-c(w)w.\nonumber
\end{eqnarray}
here $P(x,t,w,D_{x}w)=a(w,x,t)u_{xx}+b(w,x,t)u_{x}+c(w,x,t)u+h(w,x,t)$. Notice that $w\rightarrow u_{0}$ as $t$ tends to $0$ from above, implies $\mathfrak{L}_{1}(w)\rightarrow 0$ as $t\rightarrow 0^{+}$.\\

We can take $\epsilon$ small enough such that $w(x,t)$ be in a neighborhood of $u_{0}$ for all $t$ in $[0,\epsilon)$. Solving each linear problem, we end with a sequence $u^{i}$ such that $\mathcal{L}(u^{i})=(u_{0},\mathcal{L}_{1}(w^{i}))$. Consider $\mathcal{L}(u)=(u_{0},0)$ by continuity. Invoking the \emph{Inverse Function Theorem} for Banach spaces~\cite[Sec. 4.13]{zeidlerFPT}, $\mathcal{L}$ is a local diffeomorphism. Therefore, the map $\mathcal{L}(u)=(u_{0},0)$ is locally invertible. 
\end{proof}

\section{Anisotropic MCF Proposition}
This section is focused on Lagrangian methods. These methods track the evolving curve by explicitly computing the update of coordinates for each point in the curve.\\

A curve $\gamma$ evolving by MCF will stop its evolution on some point, if the curvature $\kappa$ at that point equals zero. If we need to evolve the curve further even though $\kappa=0$, we need to modify the Equation~\ref{eq:MCF} using another curvature dependent normal flow.     The equation that will describe this new flow is the following:
\begin{gather*}
\frac{\partial \gamma}{\partial t}=v\mathbf{n},
\end{gather*}
here the speed $v$ is assumed to depend on the curvature.

Consider the system:\\

Let $\rho $ be a given  distribution, and $\gamma_{0}$ be an initial smooth simple closed curve. We are interested in the evolution of $\gamma_{0}$ by\\
\begin{minipage}{0.5\linewidth}
\begin{equation*}
\begin{split}
\mbox{{\bfseries (P1)}}
  \begin{cases}\;\; \frac{\partial \gamma}{\partial t}=\frac{\partial u}{\partial \hat{\mathbf{n}}}\hat{\mathbf{n}}\\
  \\
\;\;\gamma (\cdot,0)=\gamma_{0}
\end{cases}
\end{split}
\end{equation*}
\end{minipage}
\	\
\hfill\begin{minipage}{0.5\linewidth}
\begin{equation*}
\begin{split}
\mbox{{\bfseries (P2)}}
  \begin{cases}\;\; \Delta u=\rho \;\;\;\;\; \mbox{in  }\; \Omega \\
  \\
\;\; u(\mathbf{x},0)=g(\mathbf{x})  \;\;\;\;\; \mbox{in  }\; \partial \Omega 
\end{cases}
\end{split}
\end{equation*}
\end{minipage}\\

Notice that in {\bfseries (P1)} we have a normal field which drives the evolution, and in {\bfseries (P2)}, a Poisson equation which defines a field in terms of its potential $u$ with Dirichlet's boundary condition. To link these problems with the MCF, we impose that if $\rho$ equals zero, then the MCF is recovered i.e. $\partial u/\partial{\mathbf{n}}=\kappa$. Next, we need to prove the existence and uniqueness of solutions for {\bfseries (P1)} and {\bfseries (P2)} before we introduce our numerical solution.\\

The following theorem is our main theoretical result in this paper:

\begin{theorem}\label{teo:P1P2}
Let $\gamma_{0}$ be the initial curve, denote by $\Omega$ its interior. Let $r>0$ and $\mathbf{p}$ be a point in $\Omega$ such that $|\mathbf{x}-\mathbf{p}|>0$ for all points $\mathbf{x}$ in a $r-$tubular neighborhood of $\gamma_{0}$. Let $g$ be a $C^{2}(\partial \Omega)$ function, $\rho(\mathbf{y})=\delta(\mathbf{y}-\mathbf{p}) $ be the Dirac's delta distribution, and $\gamma_{0}$ be an initial smooth simple closed curve. Then, the system ({\bfseries P1})-({\bfseries P2}) has a unique solution.
\begin{proof}
The problem ({\bfseries P2}) is a Poisson equation with a Dirichlet boundary condition. Therefore, the existence and uniqueness of its solution is well known when $g$ is a $C^{2}$ function with compact support (see, for example, ~\cite{evans2010partial}).\\

The problem ({\bfseries P1}) represents a curve being deformed in the normal direction with speed $\partial u/\partial \hat{\mathbf{n}}$. Therefore, we can consider a tubular neighborhood about $\gamma_{0}$, and verify the existence and uniqueness of solutions of ({\bfseries P1}) for small time. Similarly to Lemma~\ref{lemma:forCurves}, let $\hat{\mathbf{n}}_{0}$ be the unitary normal at time $t=0$, we reparametrize the deformations inside the $r$-tubular neighborhood by:
\begin{gather*}
\gamma (s,t)=\gamma_{0}(s)+f(s,t,u)\mathbf{\hat{n}}_{0}.
\end{gather*}
Taking the time derivative of the last equation, and recalling ({\bfseries P1}), we obtain
\begin{gather*}
\frac{\partial \gamma}{\partial t} = \left( \frac{\partial f}{\partial t}\right)\hat{\mathbf{n}}_{0}=\frac{\partial u}{\partial \mathbf{n}}\hat{\mathbf{n}}.
\end{gather*}
As $\hat{\mathbf{n}}_{0}$ is a unitary vector, we can state the partial equation for $f$:
\begin{gather}\label{ec:P1P2PDE}
\frac{\partial f}{\partial t}=\frac{\partial u}{\partial \mathbf{ n}}\langle \hat{\mathbf{n}},\hat{\mathbf{n}}_{0}\rangle.
\end{gather}
The dot product in equation (\ref{ec:P1P2PDE}) is already computed in Lemma (\ref{lemma:forCurves}), and $u$ is the solution of the Poisson problem ({\bfseries P2}). Then, invoking the general solution for Poisson's problem and the dot product in Lemma (\ref{lemma:forCurves}), we rewrite equation (\ref{ec:P1P2PDE}):
\begin{gather}\label{ec:Green}
\frac{\partial f}{\partial t}=\frac{\partial}{\partial \mathbf{n}}\left( \int_{ \Omega} G(\mathbf{x}-\mathbf{y})\,\rho(\mathbf{y})\,\mbox{d}\mathbf{\sigma_{y}}+ \int_{\partial \Omega}\frac{\partial G}{\partial n_{y}}(\mathbf{x}-\mathbf{y})g(\mathbf{y})\,d\mathbf{y}\right)\; \frac{1-\kappa_{0}}{\sqrt{(1-\kappa_{0}f)^{2}+(f^{\prime})^{2}}}
\end{gather}
The integrands in previous equation are Green's function and the Poisson kernel, whose existence is known ~\cite[Sec. 8.2]{krantz}. Comparing~\ref{ec:Green} and~\ref{ec:dfdt}, we notice that the second derivative of $f$ does not appear explicitly. Consequently, some assumptions are required in order to link ({\bfseries P1})  and ({\bfseries P2}) with the mean curvature flow: We claim that
\begin{gather}\label{ec:condition}
\frac{\partial }{\partial \mathbf{n}}\int_{\partial \Omega}\frac{\partial G}{\partial \mathbf{n}}(\mathbf{x}-\mathbf{y})g(\mathbf{y})\,d\mathbf{y}=\kappa(\mathbf{x}, t)
\end{gather}
holds, we will carry out the demonstration at the end of this proof. Substituting Equations (\ref{ec:kurvt}) and (\ref{ec:condition}) into (\ref{ec:Green}), we obtain
\begin{eqnarray}\label{ec:twocontributions}
\frac{\partial f}{\partial t}&=&\left( \frac{\partial}{\partial \mathbf{n}} \int_{ \Omega} G(\mathbf{x}-\mathbf{y})\,\rho(\mathbf{y})\,\mbox{d}\mathbf{\sigma_{y}}\right)\; \frac{1-\kappa_{0}}{\sqrt{(1-\kappa_{0}f)^{2}+(f^{\prime})^{2}}}\label{ec:evolFGreen}\\
& &+\kappa(t)\, \frac{1-\kappa_{0}}{\sqrt{(1-\kappa_{0}f)^{2}+(f^{\prime})^{2}}}\nonumber.
\end{eqnarray}
The term involving $\kappa (t)$ carries the evolution by mean curvature (compare with lemma~\ref{lemma:forCurves}). The additional term depends on the distribution function $\rho$, and does not include second order derivatives of $f$. Thus, the parabolicity is not affected. The integral in this term equals
\begin{gather*}
\int_{\Omega}G(\mathbf{x}-\mathbf{y})\rho(\mathbf{y})\,\mbox{d}\sigma_{\mathbf{y}}=\int_{\Omega}G(\mathbf{x}-\mathbf{y})\delta(\mathbf{y}-\mathbf{p})\,\mbox{d}\sigma_{\mathbf{y}}=G(\mathbf{x}-\mathbf{p}).
\end{gather*}
Recall the following properties of Green's function~\cite{evans2010partial}:
\begin{itemize}
\item $G(\mathbf{x}-\mathbf{y})$ equals zero at the boundary $\partial \Omega$.
\item Green's function satisfies -$\Delta G(\mathbf{x}-\mathbf{y})=\delta(\mathbf{x})$.
\item In general, Green's function takes the form $G(\mathbf{x}-\mathbf{y})=\Phi(|\mathbf{x}-\mathbf{y}|)+\phi(x,y)$. Here, $\Phi(|\mathbf{x}-\mathbf{y}|)$ is the \emph{fundamental solution of Laplace's equation}, $\phi(\mathbf{x},\mathbf{y})$ is known as \emph{the corrector function} and satisfies
\begin{gather*}
\Delta_{\mathbf{x}} \phi(\mathbf{x},\mathbf{y})=0 \;\; \mbox{in}\;\; \Omega\\
\phi(\mathbf{x},\mathbf{y})=\Phi(|\mathbf{x}-\mathbf{y}|)\;\; \mbox{in}\;\; \partial\Omega.
\end{gather*}
\item For a 2-D problem, 
\begin{gather*}
\Phi(|\mathbf{x}-\mathbf{y}|)=-\frac{1}{2\pi}\log(|\mathbf{x}-\mathbf{y}|).
\end{gather*}
\end{itemize}
Then, the derivatives of $G(\mathbf{x}-\mathbf{y})$ exist because $\mathbf{p}$ is outside the tubular neighborhood, and it is linearizable in this neighborhood. Consequently, equation (\ref{ec:evolFGreen}) can be linearized as a parabolic partial differential equation. The solution for the linearized problem exists by Theorem (\ref{teo:linear}), and it can be extended to a solution of (\ref{ec:evolFGreen}) by Theorem (\ref{teo:exist}).\\

Finally, we state the proof of Equation (\ref{ec:condition}), requiring that $\partial u/\partial \mathbf{n}=\kappa$ if $\rho=0$:\\

When $\rho=0$ in ({\bfseries P2}), we want to recover the MCF. Then, from problem ({\bfseries P1}) and $\rho=0$ in ({\bfseries P2}), we can state the Laplace's problem with Neumann's boundary condition
\begin{gather*}
\mbox{({\bfseries P3})}\begin{cases}\;\Delta u=0 \;\; \mbox{in}\; \Omega\\
\;\frac{\partial u}{\partial \mathbf{n}}=\kappa \; \;\mbox{in}\;\partial \Omega \;\;.
\end{cases}
\end{gather*}

Recalling problem ({\bfseries P2}), let $\mathbf{x}$ be a point in the boundary $\partial \Omega$. By Green's function relation to the Laplacian as recalled above, we write  
\begin{eqnarray}\nonumber
g(\mathbf{x})&=&u|_{\mathbf{x}\in\partial \Omega}=\int_{ \Omega}\Delta_{\mathbf{y}}G(\mathbf{x}-\mathbf{y})u(\mathbf{y})\,\mbox{d}\mathbf{y}\\
&=&\int_{\Omega}G(\mathbf{x}-\mathbf{y})\Delta u(\mathbf{y})\,\mbox{d}\mathbf{y}\;+\, \int_{\partial \Omega}\left( u(\mathbf{y})\frac{\partial G}{\partial \mathbf{n}}+G(\mathbf{x}-\mathbf{y})\frac{\partial u}{\partial \mathbf{n}}\right)\,\mbox{d}\mathbf{y}.\nonumber
\end{eqnarray}

The last equality follows form Green's identity. Substituting problem ({\bfseries P3}) yields:
\begin{eqnarray}
g(\mathbf{x})&=&u|_{\mathbf{x}\in\partial \Omega}=\int_{\partial\Omega}\left( u(\mathbf{y})\frac{\partial G}{\partial \mathbf{n}}+G(\mathbf{x}-\mathbf{{y}})\kappa(\mathbf{y})\right)\,\mbox{d}\mathbf{y}\nonumber
\end{eqnarray}
Since $G(\mathbf{x}-\mathbf{y})$ is null on the boundary, the second term cancels. Taking the normal derivative
\begin{eqnarray}\nonumber
\frac{\partial g}{\partial \mathbf{n}}&=&\frac{\partial u}{\partial \mathbf{n}}=\frac{\partial }{\partial \mathbf{n}}\int_{\partial\Omega}u(\mathbf{y})\frac{\partial G}{\partial \mathbf{n}_{\mathbf{y}}}\,\mbox{d}\mathbf{y}.
\end{eqnarray}
Using von Neumann's boundary condition in ({\bfseries P3}), one obtains (\ref{ec:condition}).\\

The uniqueness of solutions for Poisson problems determines the uniqueness of solutions for {\bfseries P2} and {\bfseries P3}. Given $\gamma_{1}(t)=(x_{1}(t),y_{1}(t))$ and $\gamma_{2}(t)=(x_{2}(t),y_{2}(t))$ two solutions of {\bfseries P1}, consider $\gamma=\gamma_{1}-\gamma_{2}$. Since  both solutions satisfy the initial data, then $\gamma(\cdot,0)=(x(0),y(0))=(0,0)$ and $\partial u/\partial \mathbf{n}=0$. This implies that $\gamma$ is stationary and $x_{1}(t)=x_{2}(t)$ and $y_{1}(t)=y_{2}(t)$.
\end{proof}
\end{theorem}

{\bfseries Remark. } A curve evolution by ({\bfseries P1}-{\bfseries P2}) subject to Equation (\ref{ec:condition}) when $\rho(\mathbf{y})=0$ states a coupled system defined by problems ({\bfseries P1}, {\bfseries P2} and {\bfseries P3}).\\

The next section deals with the numerical solution for our curve evolution. To achieve the computations, we use the integral representation for problem ({\bfseries P2}) and ({\bfseries P3}). Recalling the theory of harmonic functions and Poisson type problems~\cite{evans2010partial}, every solution of a Dirichlet problem
\begin{gather*}
\begin{cases}\;\;\Delta u=h\;\;\mbox{in}\;\;\Omega\\
u(\mathbf{x})=g(\mathbf{x})\;\;\mbox{in}\;\; \partial \Omega
\end{cases}
\end{gather*}
has a unique solution, and it must moreover satisfy \emph{Poisson's Representation Formula}:

\begin{lemma}
Let $x$ be in $\Omega$ and let
\begin{gather*}
\Phi (x) := \frac{1}{2\pi}\log \left(\frac{1}{|\mathbf{x}|}\right)
\end{gather*}
be the fundamental solution of Laplace's equation and $u$ be a solution of Poisson's equation, then the following formula holds:
\begin{gather}\label{eq:rep}
\displaystyle \int_{\Omega} \Phi(\mathbf{x}-\mathbf{y})\Delta u(\mathbf{y})d\mathbf{y}= \int_{\partial \Omega}\Phi (\mathbf{x}-\mathbf{y})\frac{ \partial u}{\partial \mathbf{n}}d\mathbf{y} -\int_{\partial \Omega}u\frac{\partial \Phi}{\partial \mathbf{n}}d\mathbf{y}-u(\mathbf{x})
\end{gather}
\end{lemma}
See~\cite{evans2010partial} for a proof.\\

In the next section, we will use a similar formula for $\mathbf{x}$ in $\partial \Omega$. These results will be crucial for our numerical approximation.

\begin{proposition}\label{prop:repre-formula}
Let $\mathbf{x}$ be in $\partial \Omega$, and $u$ be a solution of Poisson's equation, then  the following holds:
\begin{eqnarray}\label{eq:rep-f}
\int_{\Omega } \Phi(\mathbf{x}-\mathbf{y})\, \Delta u(\mathbf{y})\,d\mathbf{y}&=&\int_{\partial \Omega }\left( \Phi(\mathbf{x}-\mathbf{y})\,\frac{\partial u}{\partial \mathbf{\hat{n}}}-u\,\frac{\partial \Phi}{\partial \mathbf{\hat{n}}}(\mathbf{x}-\mathbf{y})\right)\,d\mathbf{y}-\frac{1}{2}u(\mathbf{x})\label{ec:repPoisson}
\end{eqnarray}
\end{proposition}
\begin{proof}
Notice that in Equation (\ref{eq:rep}) $\Phi$ becomes singular when $\mathbf{x}$ tends to the boundary $\partial \Omega$. Consider a neighborhood $\ell(\delta,\mathbf{x})$ over $\partial \Omega$ and then take the limit $\delta \rightarrow 0$. For the integrals on the right hand side over $\partial \Omega$ in Equation (\ref{eq:rep}) we use Green's second identity:
\begin{eqnarray}
\int_{\partial \Omega}\Phi(\mathbf{x}-\mathbf{y})\frac{\partial u}{\partial \mathbf{\hat{n}}}\,d\mathbf{y} -\int_{\partial \Omega}u\frac{\partial \Phi}{\partial \mathbf{\hat{n}}}\,d\mathbf{y}-u(\mathbf{x})&=& \int_{\partial \Omega -\ell(\delta,\mathbf{x})}\left( \Phi(\mathbf{x}-\mathbf{y})\,\frac{\partial u}{\partial \mathbf{\hat{n}}}-u\,\frac{\partial \Phi}{\partial \mathbf{\hat{n}}}(\mathbf{x}-\mathbf{y})\right)\,d\mathbf{y}\nonumber\\
&\, &\;\;\;\;\; +\cancel{\,\int_{\Omega }u\,\Delta \Phi(\mathbf{x}-\mathbf{y})\,d\mathbf{y}}\nonumber\\
&\,&\;\;\;\;\; +\int_{\ell(\delta,\mathbf{x})}\left( \Phi(\mathbf{x}-\mathbf{y})\,\frac{\partial u}{\partial \mathbf{\hat{n}}}-u\,\frac{\partial \Phi}{\partial \mathbf{\hat{n}}}(\mathbf{x}-\mathbf{y})\right)\,d\mathbf{y}.\nonumber
\end{eqnarray}
Then,

\begin{equation}\label{eq:rewrite}
\begin{aligned}
\left. \hspace{3.5em}\right.\int_{\partial \Omega}\Phi(\mathbf{x}-\mathbf{y})\frac{\partial u}{\partial \mathbf{\hat{n}}}\,d\mathbf{y} -\int_{\partial \Omega}u\frac{\partial \Phi}{\partial \mathbf{\hat{n}}}\,d\mathbf{y}-u(\mathbf{x}) = &\int\limits_{\substack{\partial \Omega -\ell(\delta,\mathbf{x})}}\left( \Phi(\mathbf{x}-\mathbf{y})\,\frac{\partial u}{\partial \mathbf{\hat{n}}}-u\,\frac{\partial \Phi}{\partial \mathbf{\hat{n}}}(\mathbf{x}-\mathbf{y})\right)\,d\mathbf{y}\\
& +\int_{\ell(\delta,\mathbf{x})} \Phi(\mathbf{x}-\mathbf{y})\,\frac{\partial u}{\partial \mathbf{\hat{n}}}\,d\mathbf{y}\\
& -\int_{\ell(\delta,\mathbf{x})}u\,\frac{\partial \Phi}{\partial \mathbf{\hat{n}}}(\mathbf{x}-\mathbf{y})\,d\mathbf{y}
\end{aligned}
\end{equation}

For the first integral over $\ell(\delta,\mathbf{x})$ in the last equation, the following limit holds:\\
\begin{eqnarray}
\lim_{\delta\rightarrow 0}\;  \int_{\ell(\delta,\mathbf{x})}\Phi(\mathbf{x}-\mathbf{y})\,\frac{\partial u}{\partial \mathbf{\hat{n}}}\,d\mathbf{y} &= & 0\nonumber
\end{eqnarray}
\indent For the second integral over $\ell(\delta,\mathbf{x})$:
\begin{eqnarray}
-\int_{\ell(\delta ,\mathbf{x})}u(\mathbf{y})\,\frac{\partial \Phi}{\partial \mathbf{\hat{n}}}(\mathbf{x}-\mathbf{y})\,d\mathbf{y}&=&-\frac{1}{2\pi\,\delta}\int_{\ell(\delta ,\mathbf{x})}u(\mathbf{y})\,d\mathbf{y}\nonumber
\end{eqnarray}
As $\delta $ tends to $0$, the last expression equals $-\frac{1}{2}u(\mathbf{x})$ in the limit. Finally, we rewrite equation (\ref{eq:rewrite}):\\
\begin{eqnarray}\label{eq:repre_f}
\int_{\Omega } \Phi(\mathbf{x}-\mathbf{y})\, \Delta u(\mathbf{y})\,d\mathbf{y}&=&\int_{\partial \Omega }\left( \Phi(\mathbf{x}-\mathbf{y})\,\frac{\partial u}{\partial \mathbf{\hat{n}}}-u\,\frac{\partial \Phi}{\partial \mathbf{\hat{n}}}(\mathbf{x}-\mathbf{y})\right)\,d\mathbf{y}-\frac{1}{2}u(\mathbf{x}).\nonumber
\end{eqnarray}
\end{proof}

\section{Numerical anisotropic MCF}\label{sec:numerical}
In this section we develop the procedure for the numerical solution of the system ({\bfseries P1}-{\bfseries P2}-{\bfseries P3}).\\

\begin{figure}[H]
\centering
\includegraphics[scale=1]{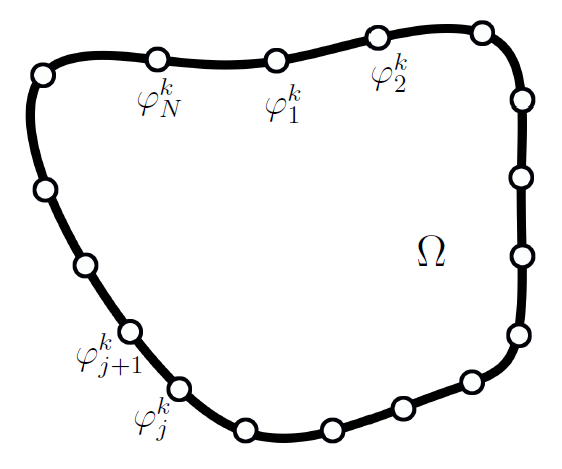}
\caption{A closed planar smooth curve is represented by a set of $N$ points $\{\varphi_{j}^{k}\}$. We also assume that the ordering index obeys the relative positions of the points.\label{fig:curve-repre}}
\end{figure}

Suppose that a curve $\varphi $ evolves by an equation of the form
\begin{gather}\label{eq:generalflow}
\frac{\partial \varphi}{\partial t}=v\mathbf{n}.
\end{gather}
This flow may be discretely approximated as follows: Let $\Delta t$ be the step in time and $k$ be a positive constant such that $k\Delta t$ represents the elapsed time. The evolving curve is represented by a set of $N$ points $\{ \varphi _{j}^{k} \}$, here $1\leq j\leq N$. For an illustrative diagram see Figure~\ref{fig:curve-repre}. The tangent vectors at $\varphi_{j}^{k}$ can be approximated using a high order finite difference formula, for example
\begin{eqnarray}\label{eq:approx1}
\left. \frac{d\varphi^{k}}{ds}\right|_{s_{j}}&\approx & \frac{-\varphi^{k}_{j+2}+8\varphi^{k}_{j+1}-8\varphi^{k}_{j-1}+\varphi^{k}_{j-2}}{12|\varphi^{k}_{j+1}-\varphi^{k}_{j}|}\nonumber\\
&=&\frac{-\varphi^{k}_{j+2}+\varphi^{k}_{j}+8\varphi^{k}_{j+1}-8\varphi^{k}_{j}-8\varphi^{k}_{j-1}+8\varphi^{k}_{j}+\varphi^{k}_{j-2}- \varphi^{k}_{j}}{12|\varphi^{k}_{j+1}-\varphi^{k}_{j}|}.
\end{eqnarray}
Using Equation (\ref{eq:approx1}) and the assumption that for all points the quantity $|\varphi_{j+1}^{k}-\varphi_{j}^{k}|$ can be approximated by a constant when fixing $k$, in~\cite{kimura} a set of formulas is presented, and they also approximate the tangent vectors $\mathbf{T}_{j}$ and curvature $\kappa_{j}$:
\begin{eqnarray}
d_{j}&:=&|\varphi_{j+1}^{k}-\varphi_{j}^{k}|\nonumber\\
\displaystyle \tau_{i}&:=&\frac{i}{|i|}\,\frac{\varphi_{j+i}^{k}-\varphi_{j}^{k}}{d_{i}}\mbox{,   for } \, i=-2,-1,1,2\nonumber\\
\displaystyle \mathbf{T}_{j}&=&\frac{-\tau_{2}+4\tau_{1}+4\tau_{-1}-\tau_{-2}}{6}\nonumber\\
\kappa_{j}&=&\mu\frac{2(\tau_{1}-\tau_{-1})}{d_{1}-d_{-1}}+(1-\mu)\frac{2(\tau_{2}-\tau_{-2})}{d_{2}-d_{-2}}\label{ec:kimuraApps}
\end{eqnarray}

Equation (\ref{eq:generalflow}) does not have a tangential term, but an approximation for tangential vectors is needed when we reparametrize the evolution in a tubular neighborhood. 
By Theorem (\ref{teo:invariance}), we know that any reparametrization will give rise to a tangential term in the evolution equation. Then, the formula for $\mathbf{T}_{j}$ provides the direction of this tangent field, we only need to know the coefficients $a_{j}^{k}$ of this field at each $ \varphi_{j}$. Notice that $|\varphi_{j+1}^{k}-\varphi_{j}^{k}|$ is a first order approximation for tangential speed. Then, rescaling $T_{j}$ formula to get $|\varphi_{j+1}^{k}-\varphi_{j}^{k}|=1$, it represents a numerical arc-length parametrization. In~\cite{kimura}, the reparametrization is such that $|\varphi_{j+1}^{k}-\varphi_{j}^{k}|$ is constant for every $j$ and fixed $k$, and the $a_{j}^{k}$ are given by
\begin{eqnarray}\label{eq:tansystem}
a_{j+1}^{k}-a_{j}^{k}&=&\frac{l_{k}/N-d_{j}^{k}}{\Delta t}\\
\displaystyle \sum_{j=1}^{n}a_{j}^{k}&=&0.\nonumber
\end{eqnarray}

When $v=\kappa$ in Equation (\ref{eq:generalflow}), that is \emph{isotropic} mean curvature flow,  an approximation for the term along the normal direction was given in~\cite{kimura}. Here we are concerned with finding a numerical solution to the system ({\bfseries P1}-{\bfseries P2}), and we seek an approximation of the form
\begin{gather}\label{solver:discretizedP1P2}
\varphi^{k+1}_{j}=\varphi_{j}^{k}+\Delta t\,\left( a_{j}^{k}\mathbf{T}_{j}^{k}+v_{j}^{k}\mathbf{n}_{j}^{k}\right).
\end{gather}

Equation (\ref{eq:tansystem}) can be used to compute the tangential component of Equation (\ref{solver:discretizedP1P2}) for the problem ({\bfseries P1}-{\bfseries P2}). The normal component requires in addition the use of Proposition (\ref{prop:repre-formula}). Substituting the data of problem ({\bfseries P2}) into equation (\ref{ec:repPoisson}), we obtain
\begin{eqnarray}\label{eq:frepre}
-\int_{\partial \Omega}\Phi(\mathbf{x}-\mathbf{y})\,\frac{\partial u(\mathbf{y})}{\partial \mathbf{n}}\,d\mathbf{y}&=&-\int_{\Omega}\Phi(\mathbf{x}-\mathbf{y})\rho(\mathbf{y})\, d\sigma_{\mathbf{y}}+\,\int_{\partial \Omega}u(\mathbf{y})\frac{\partial \Phi}{\partial \mathbf{\hat{n}}}(\mathbf{x}-\mathbf{y})\,d\mathbf{y}+\frac{1}{2}u(\mathbf{x})\;.
\end{eqnarray}
The last two terms in the right side of~\ref{eq:frepre} depend on the values of $u$ on the boundary $\partial \Omega$. According to Theorem~\ref{teo:P1P2}, we define $u$ implicitly through ({\bfseries P3}). Then, in addition to (\ref{eq:frepre}), we need to consider Poisson's representation formula for ({\bfseries P3}):
\begin{gather}
\label{ec:PoissonP3Int}
-\int_{\partial \Omega}\Phi(\mathbf{x}-\mathbf{y})\, \kappa(\mathbf{y})\,\mbox{d}\mathbf{y}=\int_{\partial \Omega}u(\mathbf{y})\frac{\partial \Phi}{\partial \mathbf{n}}(\mathbf{x}-\mathbf{y})\,\mbox{d}\mathbf{y}+\frac{1}{2}u(\mathbf{x}),
\end{gather}
Discretizing (\ref{eq:frepre}) and (\ref{ec:PoissonP3Int}) we find two systems of $N$ linear equations. Notice that, in this case, we can solve (\ref{ec:PoissonP3Int}) independently because there is no $\partial u/\partial \mathbf{n}$ dependence.\\

In order to solve equations (\ref{eq:frepre}) and (\ref{ec:PoissonP3Int}) numerically, we need to approximate the integrands on certain points. To approximate $\kappa$ at each $\varphi_{j}^{k}$, we use (\ref{ec:kimuraApps}) from Kimura's numerical scheme. In general, from discrete values at $\varphi_{j}^{k}$ and $\varphi_{j+1}^{k}$, we can determine the values for points in between by linear interpolation. This procedure is required to compute $\kappa$, $u$ and the normal vector for points along each arc $\varphi_{j},\varphi_{j+1}$. In contrast, Laplace's fundamental solution $\Phi$ and its derivatives can be substituted explicitly.\\

One can split the integrals in~(\ref{eq:frepre}) and~\ref{ec:PoissonP3Int} into $N$ integrals over the arcs $\varphi_{j},\varphi_{j+1}$. If each arc is reparametrized by an interpolating function, all the integrals can be reduced to a linear system of $N$ equations, which can be solved numerically. Then, the procedure for (\ref{eq:frepre}) and (\ref{ec:PoissonP3Int}) is to solve (\ref{ec:PoissonP3Int}) for $u$, and use these values to interpolate $u$ at the integrands in (\ref{eq:frepre}).\\

\section{Implementation of our numerical anisotropic MCF  to Contour Parametrization\label{sec:Implementation}}
For a contour parametrization problem, let $\varphi$ be a closed planar and smooth curve drawn on a picture. In addition, suppose that the curve is enclosing a unique object in the image. We perform a curve evolution that moves the curve towards the boundary of the object. We now apply the numerical procedure in the last section to perform this evolution. First, we will describe the procedure and some assumptions we will need, and finally we will state the algorithm.\\

To achieve all the numerical computations we need to define a suitable function $\rho$ for problem ({\bfseries P2}). As we pointed out previously, $\rho$ encodes the normal field which allows the curve to evolve further when $\kappa=0$. In analogy with many physical problems, suppose a point $\mathbf{p}$ is inside $\partial \Omega$, and let $\delta(\mathbf{x}-\mathbf{p})$ be the Dirac delta function. Then, let $\rho(\mathbf{y})=-\delta(\mathbf{y}-\mathbf{p})$ be the distribution in ({\bfseries P2}) and in equation~(\ref{eq:frepre}), the integral $\int_{\Omega}\Phi(\mathbf{x}-\mathbf{y})\rho(\mathbf{y})\, d\mathbf{y}$ is easily computed and curve will shrink toward $\mathbf{p}$. At each time $k\Delta t$ we need to recompute the coefficients $a_{j}^{k}$ and $v_{j}^{k}$. Thus, there is no restriction in choosing a non fixed $\mathbf{p}$, as long as it remains inside the curve $\varphi$ because the problem ({\bfseries P2}) requires a $\rho$ defined in the area $\Omega$ enclosed by $\varphi$. Although a non fixed $\mathbf{p}$ may match the curve and the object faster, a stability result in such a context would be more complicated to state. The stability of our scheme will be analyzed in the next section.\\

As the curve evolves by equation (\ref{solver:discretizedP1P2}), at each $\varphi_{j}^{k}$ we can pick a pixel value. Suppose, for illustration purposes, that the picture on which the curve is evolving is in high contrast. That is, the pixels which belong to the object have a very different value from those belonging to the background. Then, we can distinguish if a point reaches the object by its pixel value.\\

Recalling (\ref{ec:twocontributions}), the the curve evolves by a combination of the contributions of $\rho$ and a standard MCF. When a point $\varphi $ meets the object, we need to constrain its motion vanishing $\rho$ and its curvature $\kappa$. For our examples we considered black pixel objects in a white background. Let $\mbox{Pix}(\mathbf{y})$ be the pixel value function, we can make $\rho$ and $\kappa$ vanish through:
\begin{equation}\label{ec:PixConstraints}
\rho^{\ast}(\mathbf{y})=-\delta(\mathbf{y}-\mathbf{p})\frac{\mbox{Pix}(\mathbf{y})}{255},\qquad\qquad\qquad\qquad \,
\kappa^{\ast}(\mathbf{y})=\kappa(\mathbf{y})\frac{\mbox{Pix}(\mathbf{y})}{255}.
\end{equation}
Division by 255 in the previous system of equations normalizes the $\mbox{Pix}$ function, because the color intensities are usually represented by values in  $[0,255]$. Additionally, we can consider other pixel value functions depending on the number of color channels (RGB or CMYK).

To include these constraints, we summarize the previous ideas as follows:
\begin{itemize}
\item[I] Let $k=0$, $\Delta t$ be positive, consider a set of $N$ points $\{\varphi_{j}^{k}\}$ representing a closed planar smooth curve, a high contrast picture with a single object enclosed by the set of points, and an initial set of $ \{(c_{i},\mathbf{p}_{i})\}$ required to define $\rho^{\ast}(\mathbf{y})=\sum_{i}c_{i}\,\delta(\mathbf{y}-\mathbf{p}_{i})$ according to (\ref{ec:PixConstraints}).
\item[II] Check if the set of points $\mathbf{y}_{i}$ are in the enclosed area of $\{\varphi_{j}^{k}\}$. If it is not, finish procedure (or optionally, compute/ask for a new set $\{(c_{i},\mathbf{y}_{i})\}$).
 \item[III] Compute the tangents, normals, curvatures and tangential coefficients $a_{j}^{k}$ using (\ref{ec:kimuraApps}) and (\ref{eq:tansystem}).
\item [IV] Compute the normal $\partial u_{j}^{k}/\partial\mathbf{n}_{j}$ coefficients solving a $2N$ linear system:
\begin{itemize}
\item[IV.a] Discretize equation (\ref{eq:frepre}) to obtain $N$ linear equations with $2N$ unknown variables $u_{j}$ and $\partial u_{j}/\partial \mathbf{n}_{j}$. These equations cannot be modified because they represent ({\bfseries P1}) and ({\bfseries P2}). Our task is to determine the normal velocities $\partial u_{j}/\partial \mathbf{n}_{j}$.
\item[IV.b] If we use (\ref{ec:PixConstraints}) to determine where $\rho$ is zero, problem ({\bfseries P3}) takes the form: 
\begin{gather*}
\Delta u=\rho^{\ast}\;\; \mbox{on}\;\Omega ,\qquad \qquad \qquad  
\frac{\partial u}{\partial \mathbf{n}}=\kappa^{\ast} \;\; \mbox{on}\;\partial \Omega.
\end{gather*}
This system switches off $\rho$ and $\kappa$ only for those points which have zero pixel value. After discretization by (\ref{eq:frepre}), this problem will lead to $N$ linear equations which can be solved numerically  along with IV.a. 
\end{itemize}
\item[V] Compute the new positions $\{\varphi_{j}^{k+1}\} $ substituting $a_{j}^{k}$ and $v_{j}^{k}:=\partial u_{j}/\partial \mathbf{n}_{j}$ into (\ref{solver:discretizedP1P2}).
\item[VI] Check how many points have reached the object using its pixel value, if most of them have reached the object (according to the set threshold), then stop.
\item[VII] Return to II.
\end{itemize}

The stopping point in {\bf VI} will depend on a pre-assigned threshold. We used values above 90\%---explained below---for our experiments. Alternatively, one could to fix a maximum number of iterations to be carried out. This is a feature of our approach that can not be avoided, as otherwise the process will continue indefinitely.

The previous steps are distilled into the following algorithm. It has a fixed maximum possible value of $C\geq 90\%$ of contour reached, as an example.

\begin{algorithm}[H]
\caption{Contour parametrization process}\label{alg:algorithm}
\begin{algorithmic}[1]
\Procedure{}{}
\BState \emph{Initialize}:
\State $\{\varphi_{j}\}^{k=0} \gets \textit{points over initial curve}$
\State $\{(c_{i},\mathbf{y}_{i})\} \gets \textit{data to define $\rho$}$
\State $\Delta t \gets \textit{fixed step size}$
\BState \emph{loop over $k$}:
\If {$\mbox{there is } (c_{i},\mathbf{y}_{i}) \mbox{ not inside the curve}$} \State $\mathbf{STOP}$
\EndIf
\State \# $\mbox{for all }\varphi_{j}^{k}\mbox{ compute: }$
\State $\mbox{Pix}_{j}^{k}\gets\textit{pixel value}$
\State $\begin{rcases} \mathbf{T}_{j}^{k}\gets\textit{tangent vector}\;\;\;\;\;\;\\  \mathbf{n}_{j}^{k}\gets\textit{normal vector}\;\;\;\;\;\;\;\\ \mathbf{K}_{j}^{k}\gets\textit{curvature vector}\;\;\; \\ \mathbf{a}_{j}^{k}\gets\textit{tangent coefficient} \end{rcases} \mbox{using (\ref{ec:kimuraApps}) and (\ref{eq:tansystem})}$
\State $u_{j}^{k},\;(\partial u_{j}/\mathbf{n}_{j})^{k}\gets\textit{normal coefficient using (\ref{eq:frepre}), (\ref{ec:PoissonP3Int}) and (\ref{ec:PixConstraints})}$
\State $\{\varphi _{j}^{k+1}\}\gets $ \textit{ update using}\mbox{~(\ref{solver:discretizedP1P2})}
\State $C\gets \mbox{number of points matching the object}$
\If {$C\geq 90\%$} \State $\mathbf{STOP}$
\EndIf
\BState \emph{end loop}
\EndProcedure
\end{algorithmic}
\end{algorithm}

\begin{figure}
  \centering
  \includegraphics[width=0.38\textwidth]{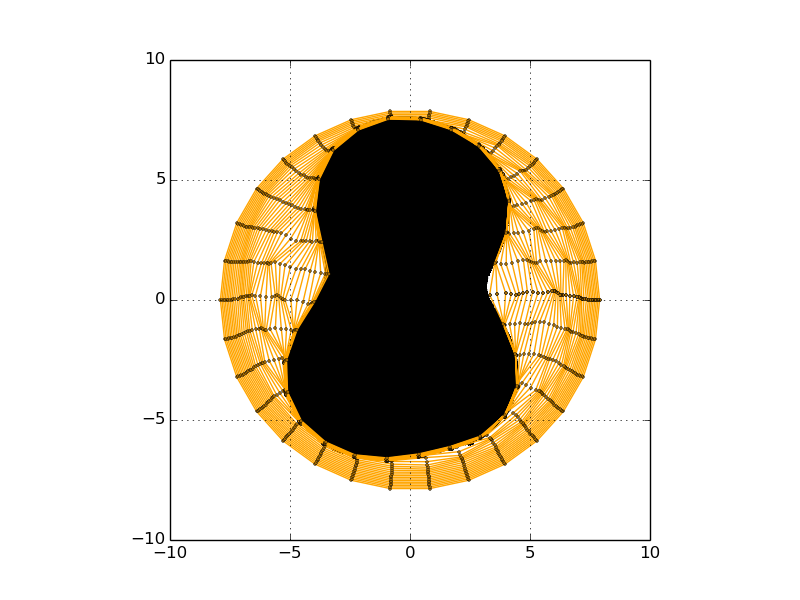}
   \includegraphics[width=0.38\textwidth]{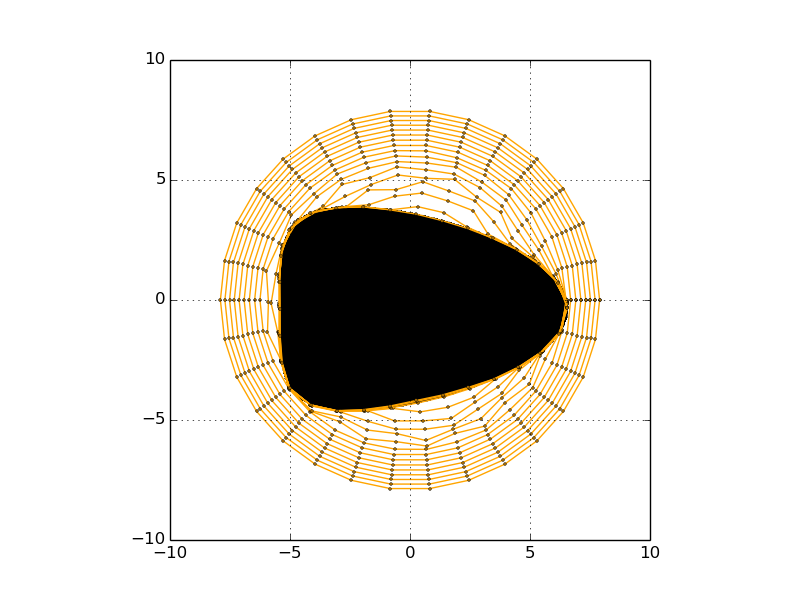}\\
   (a)\hspace{0.333\textwidth} (b)
  \caption{Non-convex and convex contour parametrization examples using our AMCF. A radius 8 circumference evolves according Algorithm~\ref{alg:algorithm} to match different pictures. In these examples, the pictures are 1800$\times$1800 pixels, scaled to 1 pixel$=$0.001, only one charge was placed at $\mathbf{p}=$(0,0). We display few curves for demonstrative purposes. The number of iterations and additional parameters are shown in Table~\ref{tab:parameters} below.\label{fig:recognition}}
  \end{figure}

We present some examples obtained with this procedure in Figure~\ref{fig:recognition} with their detailed parameters in Table~\ref{tab:parameters}. To estimate the accuracy in the contour detection, we compare the black pixels area and the area inside the last curve in evolution. The source code for these examples is available at \url{https://github.com/V3du4rd0/AMCF}.
 
\begin{table}\caption{Parameters used for pictures in Figure~\ref{fig:recognition}. Symbology: Number of points {\bfseries N}, time step {\bfseries $\Delta t$}, Parameter for curvature approximation {\bfseries $\mu$} according to~\cite{kimura}, number of iterations {\bfseries $\mathbf{\#}$it}, and the accuracy {\bfseries Acc.}($\mathbf{\%}$).\label{tab:parameters}}
\begin{threeparttable}
\begin{tabular}{lccccccc}
&  &{\bfseries N}&{\bfseries $\Delta t\times 10^{-3}$} & $\mu $&  & {\bfseries $\mathbf{\#}$ it.}&{\bfseries Acc.($\mathbf{\%}$)}\\
& & & & & &&\\[-2ex]
\cline{3-8}\\
& & & & & & &\\[-2ex]
{\small Fig.~\ref{fig:recognition}} {\bfseries a}& &30 & 1.11 & 0.15& & 26140&97\\
\multicolumn{1}{c}{{\scriptsize plotted each 1000 it. }} & & & & & &&\\
\multicolumn{1}{c}{{\scriptsize 100$\%$ of points matched}} & & & & & &&\\
& & & &  &\\\\[0.9ex]
\cline{1-8}\\
& & & &  &\\\\[0.9ex]
{\small Fig.~\ref{fig:recognition}} {\bfseries b}& &30 & 1.11  & 0.15& &30820 &98\\
\multicolumn{1}{c}{{\scriptsize plotted each 2000 it. }} & & & & & &\\
\multicolumn{1}{c}{{\scriptsize 100$\%$ of points matched}} & & & & & &&\\
& & & & & \\\\[0.9ex]
\end{tabular}
\end{threeparttable}
\end{table}

\section*{Execution times}
We used a small number of points to represent the evolving curve. Nevertheless, our numerical examples demand a considerable number of iterations. These iterations increase the computation time in a single processor computer. In order to reduce the execution time, we parallelized step 11, 12 and 13 in Algorithm~\ref{alg:algorithm}. This will reduce execution times if a long number of points are needed. All the numerical experiments were performed using double precision float point arithmetic. We used Nvidia's Jetson computer model TK1, with 2GB RAM memory, ARM A15 processor and Nvidia's Tegra K1 graphic card with 192 CUDA cores.\\

\begin{figure}
\centering
\includegraphics[scale=0.4]{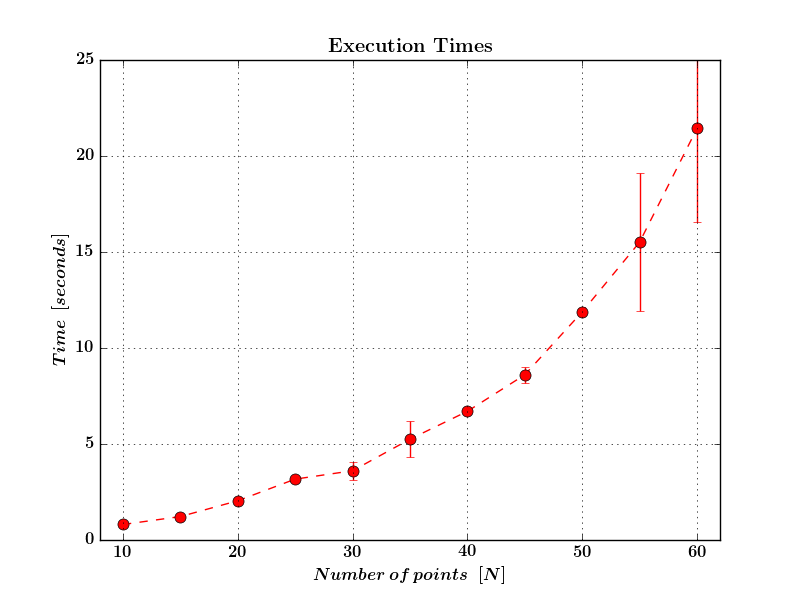}
\caption{Execution time for numerical implementation of MCF and image recognition.  \label{fig:exec}}
\end{figure}

To provide a picture of how fast our process can be performed, we implement the Algorithm (\ref{alg:algorithm}) varying the number of points $N$ and setting the time step $\Delta t=1/N^{2}$, we plot the time versus the number of points. This (execution time) plot is presented in Figure (\ref{fig:exec}). Here, the dots represent the mean time computed by taking ten samples, and the bars are the deviation of these values. The plot was obtained from a circumference as initial curve, enclosing a black circle picture with the same center. The Algorithm (\ref{alg:algorithm}) was initialized with only one $\mathbf{p}$ also localized in the center of the curve. The evolution  was stopped when at least $90\%$ of the points had reached the object. These results were achieved in parallel to speed up the process. 

\section{Conditioning and Stability of our Numerical scheme}
In this section we state our main results on the numerical analysis of the evolution scheme.\\
Let $A$ be the resulting $2N\times 2N$ matrix associated to the linear system when Equations (\ref{eq:frepre}) and (\ref{ec:PoissonP3Int}) are discretized.\\

\begin{figure}[H]
\centering
\includegraphics[height=0.35\textwidth]{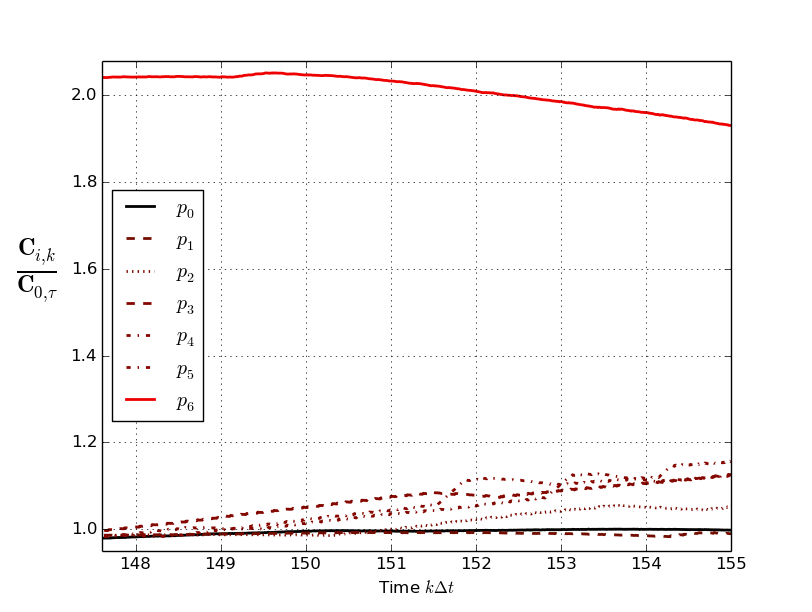}
\includegraphics[height=0.35\textwidth]{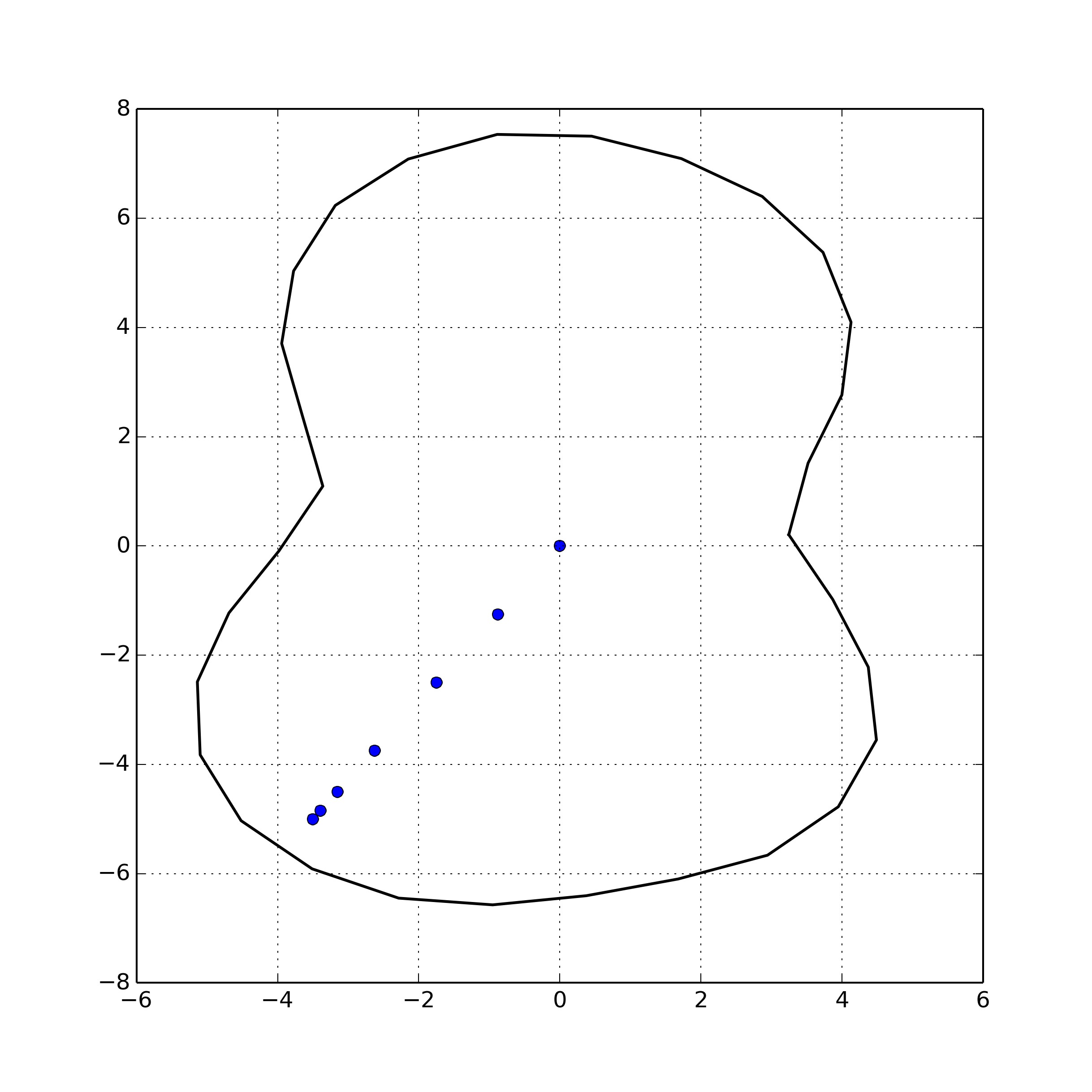}\\
(a)\hspace{0.3\textwidth}(b)
\caption{The example in figure (\ref{fig:recognition}{\bfseries a}) is repeated for different coordinates $\mathbf{p_{i}}$ from $\mathbf{p}_{0}$=(0,0) to $\mathbf{p}_{6}=$(-3.5,-5). The condition number $C(A)$ of the matrix increases as $\mathbf{p_{i}}$ approaches the curve points (a). The $x$ axis is the time when the curve matches the black pixels nearest to (-3.5,-5). In (b) we display schematically the chosen coordinates for $\mathbf{p_{i}}$.\label{fig:conditionPlot}}
\end{figure}

In order to understand the effect of the source $\mathbf{p}$, observe that equations (\ref{eq:frepre}) and (\ref{ec:PoissonP3Int}) depend on $\mathbf{p}$ through $\Phi(|\mathbf{x}-\mathbf{p}|)$. As Laplace's fundamental solution $\Phi$ is a logarithmic function, when $\mathbf{p}$ approaches the boundary, larger values of $\Phi$ and its derivatives will appear for the curve points near $\mathbf{p}$. Consequently, these values may increase some entries in the matrix in the $2N$ linear system. Therefore, the matrix norm $|| \cdot ||_{\infty}$ may change drastically when $|\mathbf{x}-\mathbf{p}|$ is lower than one because
\begin{gather*}
||A||_{\infty}=\underset{0 \leq i \leq n-1}{\max}\sum_{j=0}^{n-1}|a_{ij}|.
\end{gather*}
The condition number $C$ of a matrix $A$ with the norm $||\cdot||_{\infty}$ is thereby affected:
\begin{gather*}
C(A)=||A||_{\infty}\, ||A^{-1}||_{\infty}.
\end{gather*}
Let $\{\mathbf{p}_{i}\}_{i=0}^{n}$ be a set of source coordinates. Denote by $C_{i,k}$ the condition number of $A$ at time $k\Delta t$ when $\mathbf{p}=\mathbf{p}_{i}$. We now repeat example {\bfseries a} in Figure~\ref{fig:recognition} for different charge positions $\mathbf{p}_{i}$, from $\mathbf{p}_{0}=(0,0)$ to $\mathbf{p}_{6}=(-3.5,-5)$. In Figure (\ref{fig:conditionPlot}) we present how the condition number $C(A)$ increases as $\mathbf{p}$ approaches the curve, the condition number when $\mathbf{p}$ is at $\mathbf{p}_{0}$ is used as reference ($C_{i,k}/C_{0,k}$). Notice in Figure~\ref{fig:conditionPlot} that the highest values of $C_{i,k}$ are obtained for $i=6$, which corresponds to the point ($\mathbf{p}_{6}$) which is closest to boundary.

{\bfseries Remark.} In order to maximize the distance between $\mathbf{p}$ and the curve in a convex contour parameterization problem, a suitable choice of $\mathbf{p}$ at time $k\Delta t$ is the centroid coordinates of the set $\{\varphi_{j}^{k} \}_{j=0}^{N-1}$. Further investigating other possible choices of placement for the charge in the general case, or increasing the number of charges, constitutes an interesting avenue for future research, and lies outside the scope of our work here.\\

We now focus on the stability of the updating scheme
\begin{gather}\label{ec:scheme}
\varphi^{k+1}_{j}=\varphi_{j}^{k}+\Delta t\,\left( a_{j}^{k}\mathbf{T}_{j}^{k}+v_{j}^{k}\mathbf{n}_{j}^{k}\right).   
\end{gather}
Let $\varphi(s,t)$ be the curve at time $t=k\,\Delta t$ represented by a set of $N$ points $\{ \varphi_{j}^{k}\}_{j=0}^{N-1}$, its length $L$ is approximated by $\ell=\sum_{j=1}^{N}|\varphi_{j+1}-\varphi_{j}|$. We will rewrite the tangential and normal vectors in~(\ref{ec:scheme}) and thus obtain a expression only in terms of the set $\{\varphi_{j}^{k}\}_{j=0}^{N-1}$. Recall the definitions of $\mathbf{T}_{j}^{k}$ and $\mathbf{n}_{j}^{k}$:
\begin{eqnarray}\label{ec:relations}
\mathbf{T}_{j}^{k}&:=&\frac{-\tau_{2}+4\tau_{1}+4\tau_{-1}-\tau_{-2}}{6}\nonumber\\
\tau_{i}&:=&\mbox{sign}(i)\frac{\varphi_{j+i}^{k}-\varphi_{j}^{k}}{d_{i}}\nonumber\\
d_{i}&:=&|\varphi_{j+i}^{k}-\varphi_{j}^{k}|\\
\mathbf{n}_{j}^{k}&:=&\left(\mathbf{T}_{j}^{k}\right)^{\bot}=\begin{pmatrix}
0 & -1 \\ 1 & 0
\end{pmatrix}\mathbf{T}_{j}^{k}\nonumber
\end{eqnarray}
Therefore
\begin{eqnarray}\label{ec:Tangential}
\mathbf{T}_{j}^{k}&=&\frac{1}{6}\left( -\frac{(\varphi_{j+2}-\varphi_{j})}{d_{2}}+\frac{4(\varphi_{j+1}-\varphi_{j})}{d_{1}}-\frac{4(\varphi_{j-1}-\varphi_{j})}{d_{-1}}+\frac{(\varphi_{j-2}-\varphi_{j})}{d_{-2}}\right).
\end{eqnarray}

To determine the numerical error propagation in~(\ref{ec:scheme}), let $\{\varphi_{\ast\,s}^{k}\}_{s=0}^{N_{1}-1}$ be an exact solution of system ({\bfseries P1}), ({\bfseries P2}) and ({\bfseries P3}), then all points $\varphi_{\ast\,s}$ must satisfy~(\ref{ec:scheme}). If another initial representation of the same curve evolves, the final positions in both sets may represent different curves, because the approximation errors in each set propagate differently. The core in the following analysis is to compare the error amplification when the relative distances $|\varphi_{s+1}^{k}-\varphi_{s}^{k}|$ change, and this is crucial for the stability as the evolving curves shrink along the flow.\\

Let $\{ \varphi_{\ast\,s}\}_{s=0}^{N_{1}-1}$ be the exact solution, and let $\{\varphi_{s}\}_{s=0}^{N_{2}-1} $ be another discrete representation of the same evolving curve with uniformly distributed points. Suppose that both sets have one point $\varphi_{j}$ in common. The difference between these sets is the number of points, this condition prevents us to compare the coordinates point-wise, except for the common point $\varphi_{j}$ and when the same time has elapsed.\\

As $\varphi_{\ast}$ is an exact solution of the system ({\bfseries P1}), ({\bfseries P2}) and ({\bfseries P3}), then it must satisfy~(\ref{ec:scheme}). At the time $(k+1 )\Delta t$, the coordinates for the common point $\varphi_{j}$ will differ by a small error $\varepsilon_{j}^{k}$ 
\begin{gather}\label{ec:difference}
\varphi_{j}^{k+1}=\varphi_{\ast\,j}^{k+1}+\varepsilon(\varphi_{j}^{k+1}),
\end{gather}

denoting $\varepsilon(\varphi_{j}^{k})$ by $\varepsilon_{j}^{k}$, we point out that $\varepsilon_{s}^{k}\neq \vec{0}$ except for $s=j$. We will get an equation for $\varepsilon_{j}^{k+1}$ in terms on the error in the additional parameters: 
\begin{gather*}
a_{j}^{k}=a_{\ast\,j}^{k}+a_{\varepsilon\,j} ^{k}\hspace{0.3\textwidth} v_{j}^{k}=v_{\ast\,j}^{k}+v_{\varepsilon\,j}^{k} \\
\mathbf{n}_{j}^{k}=\mathbf{n}_{\ast\,j}^{k}+n_{\varepsilon\,j}^{k}
\end{gather*}
and so on for $\tau_{i}$ and $\mathbf{T}_{j}^{k}$ in $\varphi_{j}^{k}$. Notice that, even when $\varphi_{j}^{k}=\varphi_{\ast\,j}^{k}$, and therefore $\varepsilon_{j}^{k}=0$, the previous parameters depend non linearly on other points. Thus, the coefficients which drive the evolution ($a_{\cdot}^{k}$ and $v_{\cdot}^{k}$) are different for $\{\varphi_{s}^{k}\}_{s=0}^{N_{2}-1}$ and $\{\varphi_{\ast\,s}^{k}\}^{N_{1}-1}$.\\

Since both sets of points are uniformly distributed at time $k\Delta t$, then $a_{j}^{k}=a_{\ast\,j}^{k}=0 $. We can now rewrite Equation~\ref{ec:scheme}:
\begin{gather*}
\varphi_{\ast\, j}^{k+1}+\varepsilon_{j}^{k+1}=\varphi_{\ast\, j}^{k}+\cancel{\varepsilon_{j}^{k}}+\Delta t \left( \cancel{(a_{\ast\,j}^{k}+a_{\varepsilon\,j})(\mathbf{T}_{\ast\,j}^{k}+T_{\varepsilon\,j})}+(v_{\ast\,j}+v_{\varepsilon})(\mathbf{n}_{\ast\,j}+n_{\varepsilon\,j})\right).
\end{gather*}
{\bfseries Remark. } We cancel the tangential terms because the points are uniformly distributed at time $k$, but in further iterations, for example  $\varphi_{j}^{k+2}$, this is not possible in general. Consequently, this analysis determines whether choosing a different number of points is appropriate to reduce the error amplification only for the next time step. If we do not assume this uniform tangential distribution, the resulting equation will be complicated to analyze, even when~(\ref{ec:scheme}) is linearized, because it would then depend on the projection of $\varepsilon_{s}$ along the segments $\varphi_{\ast\,s},\varphi_{\ast\,s+1}$ for all $s=0,\cdots ,N-1$. \\

Recall that $\varphi_{\ast\,j}^{k}$ is an exact solution, so the last equation implies:
\begin{gather}\label{ec:error}
\varepsilon_{j}^{k+1}=\Delta t \left( v_{\varepsilon\,j}^{k}\mathbf{n}_{\ast\,j}^{k}+v_{\ast\,j}^{k}n_{\varepsilon\, j}^{k} +v_{\epsilon\,j}^{k}n_{\varepsilon\,j}^{k}\right).
\end{gather}
To obtain an explicit expression for~(\ref{ec:error}), we need the following computations:
\begin{eqnarray}
T_{\varepsilon\,j}^{k}&=&\mathbf{T}_{j}^{k}-T_{\ast\,j}^{k}\nonumber\\
&=&\frac{1}{6}\left( -\frac{(\varphi_{\ast\,j+2}-\varphi_{\ast\,j})}{d_{2}}+\frac{4(\varphi_{\ast\,j+1}-\varphi_{\ast\,j})}{d_{1}}-\frac{4(\varphi_{\ast\,j-1}-\varphi_{\ast\,j})}{d_{-1}}\right)\nonumber\\
& &+\frac{1}{6}\left(\frac{(\varphi_{\ast\,j-2}-\varphi_{\ast\,j})}{d_{-2}}-\frac{\varepsilon_{j+2}}{d_{2}}+\frac{4\varepsilon_{j+1}}{d_{1}}-\frac{4\varepsilon_{j-1}}{d_{1}}+\frac{\varepsilon_{-2}}{d_{-2}}\right)\nonumber\\
& &-\frac{1}{6}\left( -\frac{(\varphi_{\ast\,j+2}-\varphi_{\ast\,j})}{d_{\ast\,2}}+\frac{4(\varphi_{\ast\,j+1}-\varphi_{\ast\,j})}{d_{\ast\,1}}-\frac{4(\varphi_{\ast\,j-1}-\varphi_{\ast\,j})}{d_{\ast\,-1}}+\frac{(\varphi_{\ast\,j-2}-\varphi_{\ast\,j})}{d_{\ast\,-2}}\right).\nonumber
\end{eqnarray}
Assuming that the difference between $N_{1}$ and $N_{2}$ is small, we linearize  $1/d_{i}$. Then,
\begin{eqnarray}
T_{\varepsilon\,j}^{k}&=&\frac{1}{6}\left( \frac{(\varphi_{\ast\,j+2}-\varphi_{\ast\,j})\cdot \varepsilon_{j+2}}{d_{\ast\,2}^{3}}(\varphi_{\ast\,j+2}-\varphi_{\ast\,j})-\frac{(\varphi_{\ast\,j+1}-\varphi_{\ast\,j})\cdot \varepsilon_{j+1}}{d_{\ast\,1}^{3}}(\varphi_{\ast\,j+1}-\varphi_{\ast\,j})\right.\nonumber\\
& &+\frac{(\varphi_{\ast\,j-1}-\varphi_{\ast\,j})\cdot \varepsilon_{j-1}}{d_{\ast\,-1}^{3}}(\varphi_{\ast\,j-1}-\varphi_{\ast\,j})-\frac{(\varphi_{\ast\,j-2}-\varphi_{\ast\,j})\cdot \varepsilon_{j-2}}{d_{\ast\,-2}^{3}}(\varphi_{\ast\,j-2}-\varphi_{\ast\,j})\nonumber\\
& &\left. -\frac{\varepsilon_{j+2}}{d_{2}}+\frac{4\varepsilon_{j+1}}{d_{1}}-\frac{4\varepsilon_{j-1}}{d_{1}}+\frac{\varepsilon_{-2}}{d_{-2}} \right).\nonumber
\end{eqnarray}

To simplify a future computation, we write the vectors $\varphi_{\ast\,j+i}-\varphi_{\ast\,j}$ in terms of $ \varepsilon_{j+i}$ and $\varepsilon_{j+i}^{\perp}$:
\begin{eqnarray}
T_{\varepsilon\,j}^{k}&=&C_{2}\,\varepsilon_{2}+C_{1}\,\varepsilon_{1}+C_{-1}\,\varepsilon_{-1}+C_{-2}\,\varepsilon_{-2}\nonumber\\
& &+D_{2}\,\varepsilon_{2}^{\perp}-D_{1}\,\varepsilon_{1}^{\perp}-D_{-1}\,\varepsilon_{-1}^{\perp}+D_{-2}\,\varepsilon_{-2}^{\perp}\nonumber,
\end{eqnarray}
here
\begin{gather*}
C_{\pm 2}=\frac{\mbox{{\bfseries sign}}(\pm 2)}{6}\left( \frac{\langle \varphi_{\ast\,j\pm 2}-\varphi_{\ast\,j}\, ,\, \varepsilon_{j\pm 2}\rangle^{2}}{|\varphi_{\ast\,j\pm 
2}-\varphi_{\ast\,j}|^{3}|\varepsilon_{j\pm 2}|^{2}}-\frac{1}{|\varphi_{j\pm 2}-\varphi_{j}|}\right)\\
C_{\pm 1}=\frac{\mbox{{\bfseries sign}}(\pm 1)}{6}\left( -\frac{\langle \varphi_{\ast\,j\pm 1}-\varphi_{\ast\,j}\, ,\, \varepsilon_{j\pm 1}\rangle^{2}}{|\varphi_{\ast\,j\pm 1}-\varphi_{\ast\,j}|^{3}|\varepsilon_{j\pm 1}|^{2}}+\frac{4}{|\varphi_{j\pm 1}-\varphi_{j}|}\right)\\
D_{\pm i}=\frac{\mbox{{\bfseries sign}}(\pm 1)}{6}\frac{\langle \varphi_{\ast\,j\pm i}-\varphi_{\ast\,j}\, ,\, \varepsilon_{j\pm i}\rangle \langle \varphi_{\ast\,j\pm i}-\varphi_{\ast\,j}\, ,\, \varepsilon_{j\pm i}^{\perp}\rangle}{|\varphi_{\ast\,j\pm i}-\varphi_{\ast\,j}|^{3}|\varepsilon_{j\pm i}|^{2}}.
\end{gather*}
Additionally, we write $\mathbf{n}_{\ast\,j}^{k}$ in terms of $\varepsilon_{j+i}$ and $\varepsilon_{j+i}^{\perp}$:
\begin{eqnarray}
\mathbf{n}_{\ast\,j}^{k}&=&A_{2}\varepsilon_{j+2}+4A_{1}\varepsilon_{j+1}+4A_{-1}\varepsilon_{j-1}+A_{-2}\varepsilon_{j-2}\nonumber\\
& &+B_{2}\varepsilon_{j+2}^{\perp}+4B_{1}\varepsilon_{j+1}^{\perp}+4B_{-1}\varepsilon_{j-1}+B_{-2}\varepsilon_{j-2}^{\perp}\nonumber
\end{eqnarray}
\begin{gather*}
A_{\pm i}=\frac{\mbox{{\bfseries sign}}(\pm i)}{6}\frac{\langle \varphi_{\ast\,j\pm i}-\varphi_{\ast\,j}\, ,\, \varepsilon_{j\pm i}\rangle}{|\varphi_{\ast\,j\pm i}-\varphi_{\ast\,j}|\,|\varepsilon_{j\pm i}|^{2}} \hspace{0.15\textwidth}B_{\pm i}=\frac{\mbox{{\bfseries sign}}(\pm i)}{6}\frac{\langle \varphi_{\ast\,j\pm i}-\varphi_{\ast\,j}\, ,\, \varepsilon_{j\pm i}^{\perp}\rangle}{|\varphi_{\ast\,j\pm i}-\varphi_{\ast\,j}|\,|\varepsilon_{j\pm i}|^{2}}
\end{gather*}
Later on, we will bound these constants. Let $\mathbf{J}=\begin{pmatrix} 0 & -1\\ 1 & 0\end{pmatrix}$ and $\mathbf{I}$ be the identity matrix, define\\
\begin{minipage}{0.5\textwidth}
\begin{gather*}
\alpha_{\pm 2}=v_{\varepsilon\,j} A_{\pm 2}+(v_{\ast\,j}+v_{\varepsilon\,j})D_{\pm 2}\\
\beta_{\pm 2}=v_{\varepsilon\,j}B_{\pm 2}+(v_{\ast\,j}+v_{\varepsilon\,j})C_{\pm 2}
\end{gather*}
\end{minipage}\begin{minipage}{0.5\textwidth}
\begin{gather*}
\alpha_{\pm 1}=4A_{\pm 1}v_{\varepsilon\, j}-\mathbf{sign}(\pm 1)\,(v_{\ast\,j}+v_{\varepsilon\,j})D_{\pm 1}\\
\beta_{\pm 1}=4B_{\pm 1}v_{\varepsilon\,j}+(v_{\ast\,j}+v_{\varepsilon\,j})C_{\pm 1},
\end{gather*}
\end{minipage}\\
we can now rewrite Equation~\ref{ec:error}
\begin{eqnarray}
\varepsilon_{j}^{k+1}&=&\Delta t \left( v_{\varepsilon\,j}^{k}\mathbf{n}_{\ast\,j}^{k}+v_{j}^{k}n_{\varepsilon\, j}^{k}+v_{\varepsilon\,j}^{k}n_{\varepsilon\, j}^{k}\right)\nonumber\\
&=&\Delta t\left( (\alpha_{ 2}\mathbf{I}+\beta_{ 2}\mathbf{J} ) \varepsilon_{j+2}+(\alpha_{ 1}\mathbf{I}+\beta_{ 1}\mathbf{J}) \varepsilon_{j+1}+(\alpha_{-1}\mathbf{I}+\beta_{-1}\mathbf{J}) \varepsilon_{j-1}+(\alpha_{-2}\mathbf{I}+\beta_{-1}\mathbf{J}) \varepsilon_{j-2}\right).\label{ec:evol-error}
\end{eqnarray}

As the curve is closed, we can assume that $\varepsilon$ is periodic of period $L$ over the curve. Von Neumann's theory implements Fourier analysis to determine whether the error $\varepsilon_{j}$ remains bounded as time elapses. The discrete representation of the curve requires us to consider the \emph{discrete Fourier transformation} of $\varepsilon$ in the spatial coordinate and with periodic conditions. From Fourier analysis, $\varepsilon$ can be written using the inverse Fourier transformation as\\
\begin{minipage}{0.5\textwidth}
\begin{gather*}
\varepsilon(\varphi_{j}^{k})=\frac{1}{N_{2}}\sum_{m=0}^{N_{2}-1}e^{i\,\frac{2\pi}{N_{2}}\,j\,m}b_{m}(k),
\end{gather*}
\end{minipage}
\begin{minipage}{0.5\textwidth}
\begin{gather*}
b_{m}(k)=\sum_{m=0}^{N_{2}-1}e^{-i\frac{2\pi}{N_{2}}\,j\,m}\varepsilon(\varphi_{j}^{k}).
\end{gather*}
\end{minipage}\\
Thus
\begin{gather}\label{ec:DFT2}
\varepsilon_{j+1}^{k}=\frac{1}{N_{2}}\sum_{m=0}^{N_{2}-1}e^{i\frac{2\pi}{N_{2}}\,jm}e^{i\frac{2\pi}{N_{2}}\,m}b_{m}(k).
\end{gather}
After substituting into~(\ref{ec:evol-error}) we find
\begin{eqnarray}\label{ec:ampl_E}
\frac{1}{N_{2}}\sum_{m=0}^{N-1}e^{i\frac{2\pi}{N_{2}}\,j \,m}b_{m}(k+1)&=&\frac{1}{N_{2}}\sum_{m=0}^{N-1}\begin{pmatrix}
\zeta & -\eta\\
\eta & \zeta
\end{pmatrix} e^{i\frac{2\pi}{N_{2}}\,j \,m}b_{m}(k).
\end{eqnarray}
Here
\begin{equation*}
\begin{aligned}
\zeta&=\Delta t\left(e^{2i\frac{2\pi}{N_{2}}m}\alpha_{2}+e^{i\frac{2\pi}{N_{2}}m}\alpha_{1}+e^{-i\frac{2\pi}{N_{2}}m}\alpha_{-1}+e^{-2i\frac{2\pi}{N_{2}}m}\alpha_{-2}\right)\\
\eta&=\Delta t\left(e^{2i\frac{2\pi}{N_{2}}m}\beta_{2}+e^{i\frac{2\pi}{N_{2}}m}\beta_{1}+e^{-i\frac{2\pi}{N_{2}}m}\beta_{-1}+e^{-2i\frac{2\pi}{N_{2}}m}\beta_{-2}\right).
\end{aligned}
\end{equation*}
Equation~(\ref{ec:ampl_E}) defines a formula from which we deduce $b_{m}(k)$ in terms of $b_{m}(k-1)$. However, this formula is not recursive because the matrix entries depend on time through $\alpha_{i}$ and $\beta_{i}$. Let $M$ be given by
\begin{gather*}
M:=\begin{pmatrix}
\zeta & -\eta\\
\eta & \zeta
\end{pmatrix}
\end{gather*}

\begin{eqnarray}\label{eq:recursive}
b_{m}(k)&=&Mb_{m}(k-1).\nonumber
\end{eqnarray}

Then, using von Neumann's analysis we can only determine how much does the error is being amplified at the next time step. To proceed with this analysis, we need for $\varepsilon$ to remain bounded in the consecutive step. Therefore, we will compute the eigenvalues $\lambda_{1}$ and $\lambda_{2}$ of $M$. Letting  $w_{1}$ and $w_{2}$ be the associated eigenvectors, then $\varepsilon$ can be written in terms of $w_{1}$ and $w_{2}$. Consequently, we will guarantee that $\varepsilon $ does not grow by finding the constraints on $N_{2}$, such that the eigenvalues are bounded absolutely by one.\\

Let $\theta:=\frac{2\pi}{N_{2}}\,m$, and let $\lambda_{1}$, $\lambda_{2}$ be the eigenvalues of $M$. Then, the eigenvalue moduli are
\begin{gather}
|\lambda_{1}|^{2}=\Delta t^{2}\left( \sum\limits_{\substack{i=-2 \\ i\neq 0}}^{2}\left(\cos \left( i\,\theta\right)\alpha_{i}-\sin \left( i\,\theta\right)\beta_{i}\right) \right)^{2}+\Delta t^{2}\left( \sum\limits_{\substack{i=-2 \\ i\neq 0}}^{2}\left(\sin \left( i\,\theta\right)\alpha_{i}+\cos \left( i\,\theta\right)\beta_{i}\right) \right)^{2}\label{ec:eigenMod1}
\end{gather}
\begin{gather}\label{ec:eigenMod2}
|\lambda_{2}|^{2}=\Delta t^{2}\left( \sum\limits_{\substack{i=-2 \\ i\neq 0}}^{2}\left(\cos \left( i\,\theta\right)\alpha_{i}+\sin \left( i\,\theta\right)\beta_{i}\right) \right)^{2}+\Delta t^{2}\left( \sum\limits_{\substack{i=-2 \\ i\neq 0}}^{2}\left(\sin \left( i\,\theta\right)\alpha_{i}-\cos \left( i\,\theta\right)\beta_{i}\right) \right)^{2}
\end{gather}

From now, we will only consider the amplification of the worse possible error. This error amplification can be obtained by considering the maximum of $|\varepsilon_{i}|$ and replacing the values of sine and cosine above by 1 or -1, in that order.\\

As the points lie equidistantly, and do not overlap we find:\\
\begin{minipage}{0.5\textwidth}
\begin{gather*}
|\varphi_{\ast\,j+1}-\varphi_{j}|=|\varphi_{\ast\,j-1}-\varphi_{j}|=\frac{L}{N_{1}}\\
\frac{L}{N_{1}}<|\varphi_{\ast\,j\pm 2}-\varphi_{j}|\leq 2\frac{L}{N_{1}}
\end{gather*}
\end{minipage}\begin{minipage}{0.5\textwidth}
\begin{gather*}
|\varphi_{j+1}-\varphi_{j}|=|\varphi_{j-1}-\varphi_{j}|=\frac{L}{N_{2}}\\
\frac{L}{N_{2}}<|\varphi_{j\pm 2}-\varphi_{j}|\leq 2\frac{L}{N_{2}}.
\end{gather*}
\end{minipage}\\

\noindent In addition,
\begin{gather*}
L\left| \frac{1}{N_{1}}-\frac{1}{N_{2}}\right|=\left|\frac{}{} \,|\varphi_{\ast\,j\pm 1}-\varphi_{\ast\,j}|-|\varphi_{j\pm 1}-\varphi_{j}|\,\right|\leq\left| \varphi_{\ast\,j\pm 1}-\cancel{\varphi_{\ast\,j}}-\varphi_{j\pm 1}+\cancel{\varphi_{j}}\right|=|\varepsilon_{j\pm 1}|.
\end{gather*}

We consider the following bounds using $\langle  \frac{\varphi_{\ast\,j+i}-\varphi_{j}}{|\varphi_{\ast\,j+i}-\varphi_{j}|},\varepsilon_{j+i}\rangle\leq |\varepsilon_{j+i}| $ and $\langle  \frac{\varphi_{\ast\,j+i}-\varphi_{j}}{|\varphi_{\ast\,j+i}-\varphi_{j}|},\varepsilon_{j+i}^{\perp}\rangle\leq |\varepsilon_{j+i}| $\\
\begin{minipage}{0.33\textwidth}
\begin{gather*}
|A_{i}|\leq\frac{1}{6\,|\varepsilon_{j+i}|}
\end{gather*}
\end{minipage}\begin{minipage}{0.33\textwidth}
\begin{gather*}
|B_{i}|\leq\frac{1}{6\,|\varepsilon_{j+i}|}
\end{gather*}
\end{minipage}\begin{minipage}{0.33\textwidth}
\begin{gather*}
|D_{i}|\leq\frac{1}{6\,|\varphi_{\ast\,j+i}-\varphi_{\ast\,j}|}
\end{gather*}
\end{minipage}\\
\begin{eqnarray}
|C_{1}|+|C_{-1}|&=&\frac{1}{6}\left( \left |-\frac{\langle \varphi_{\ast\,j+1}-\varphi_{\ast\,j},\varepsilon_{j+1}\rangle^{2}}{|\varphi_{\ast\,j+1}-\varphi_{\ast\,j}|^{3}\,|\varepsilon_{j+1}|^{2}}+\frac{4}{|\varphi_{j+1}-\varphi_{j}|}\right|   +\left|\frac{\langle \varphi_{\ast\,j-1}-\varphi_{\ast\,j},\varepsilon_{j-1}\rangle^{2}}{|\varphi_{\ast\,j-1}-\varphi_{\ast\,j}|^{3}\,|\varepsilon_{j-1}|^{2}}-\frac{4}{|\varphi_{j-1}-\varphi_{j}|}\right|\right)\nonumber\\
&\leq &\frac{1}{6}\max\left\{ \frac{1}{|\varphi_{\ast\,j+1}-\varphi_{\ast\,j}|},\frac{4}{|\varphi_{j+1}-\varphi_{j}|}\right\}+\frac{1}{6}\max\left\{ \frac{1}{|\varphi_{\ast\,j-1}-\varphi_{\ast\,j}|},\frac{4}{|\varphi_{j-1}-\varphi_{j}|}\right\}\nonumber\\
&=&\frac{1}{3L}\max\{ N_{1},4N_{2}\}\nonumber\\
|C_{2}|+|C_{-2}|&=&\frac{1}{6}\left(\left| \frac{\langle \varphi_{\ast\,j+2}-\varphi_{\ast\,j},\varepsilon_{j+2}\rangle^{2}}{|\varphi_{\ast\,j+2}-\varphi_{\ast\,j}|^{3}\,|\varepsilon_{j+2}|^{2}}-\frac{1}{|\varphi_{j+2}-\varphi_{j}|}\right|   +\left|-\frac{\langle \varphi_{\ast\,j-2}-\varphi_{\ast\,j},\varepsilon_{j-2}\rangle^{2}}{|\varphi_{\ast\,j-2}-\varphi_{\ast\,j}|^{3}\,|\varepsilon_{j-2}|^{2}}+\frac{1}{|\varphi_{j-2}-\varphi_{j}|}\right|\right)\nonumber\\
&\leq &\frac{1}{6}\max\left\{ \frac{1}{|\varphi_{\ast\,j+2}-\varphi_{\ast\,j}|},\frac{1}{|\varphi_{j+2}-\varphi_{j}|}\right\}+\frac{1}{6}\max\left\{ \frac{1}{|\varphi_{\ast\,j-2}-\varphi_{\ast\,j}|},\frac{1}{|\varphi_{j-2}-\varphi_{j}|}\right\}\nonumber\\
&=& \frac{1}{3L}\max\left\{N_{1},N_{2} \right\}\nonumber
\end{eqnarray}

Recalling~(\ref{ec:eigenMod1}) and~(\ref{ec:eigenMod2}), we now focus on a bound for $\sum\limits_{\substack{i=-2\\ i \neq 0}}^{2}\alpha_{i}+\beta_{i}$.
\begin{eqnarray}
\alpha_{1}+\alpha_{-1}+\beta_{1}+\beta_{-1}& = & 4v_{\varepsilon\, j}(A_{1}+A_{-1}+B_{1}+B_{-1})+(v_{\ast\,j}+v_{\varepsilon\,j})(C_{1}+C_{-1}+D_{1}-D_{-1})\nonumber\\
&\leq &\frac{8\,|v_{\varepsilon\,j}|}{3\,\min(|\varepsilon_{j\pm 1}|)}+\frac{|v_{\ast\,j}+v_{\varepsilon\,j}|}{3L}\left(\max\{N_{1},4N_{2}\}+N_{2}\right)\nonumber\\
\alpha_{2}+\alpha_{-2}+\beta_{2}+\beta_{-2}&= & v_{\varepsilon\, j}(A_{2}+A_{-2}+B_{2}+B_{-2})+(v_{\ast\,j}+v_{\varepsilon\,j})(C_{2}+C_{-2}+D_{2}-D_{-2})\nonumber\\
&\leq &\frac{|v_{\varepsilon\,j}|}{3\min(|\varepsilon_{j\pm 2}|)}+\frac{|v_{\ast\,j}+v_{\varepsilon\,j}|}{3L}(\max\{ N_{1},N_{2}\}+2N_{2})\nonumber
\end{eqnarray}

Assuming that the size of errors $|\varepsilon_{j\pm i}|$ is almost the same for $i=-2, -1,1,2$. Then, the norm of the eigenvalues in~(\ref{ec:eigenMod1}) and~(\ref{ec:eigenMod2}) can be bounded by
\begin{gather}
|\lambda_{\cdot}|^{2}\leq 2\Delta t^{2}\left( \frac{9\,v_{\varepsilon\,j}}{3\min(|\varepsilon_{j\pm i}|)}+\frac{|v_{\ast\,j}+v_{\varepsilon\,j}|}{3L}(2\max\{N_{1},4N_{2}\}+3N_{2})\right)^{2}.\label{ec:eigenBound}
\end{gather}
{\bfseries Remark. } The number of points that represents the curve determines the integral discretization, because we split the integral over $\partial \Omega$ into $N$ integrals over the arcs $\varphi_{s},\varphi_{s+1}$. However, all these integrals are computed by an approximation formula, and therefore this error propagates to $v_{j}$, we called this error $v_{\varepsilon}$. Some well-known accurate formulas can be used to bound this integration error, in all numerical experiments we used the 3-points \emph{Gaussian Quadrature Formula}. The accuracy of these formulas increase with the number of interpolated points. To proceed with the analysis, we will disregard $v_{\varepsilon\,j}$ in~(\ref{ec:eigenBound}).\\

\begin{proposition}\label{prop:criteria1}
Let $\{ \varphi_{\ast\,s}\}_{s=0}^{N_{1}-1}$ and $\{\varphi_{j}\}_{s=0}^{N_{2}-1}$ be two discrete representations of the same curve evolving by {\bfseries P1}, {\bfseries P2} and {\bfseries P3}. Then the representation with the least number of points will produce lower error amplification.
 \begin{proof}
Using Equation~\ref{ec:eigenBound}, and dealing with $v_{\varepsilon\,j}$ according to the last remark:
\begin{eqnarray}
|\lambda_{\cdot}|&\leq& 2\Delta t^{2}\left( \frac{|v_{\ast\,j}|}{3L}(2\max\{N_{1},4N_{2}\}+3N_{2})\right)^{2}.\nonumber
\end{eqnarray}

We are using Kimura's uniform tangential redistribution scheme, so we are constrained to set $\Delta t=1/N_{1}^{2}$. Therefore,
\begin{eqnarray}
|\lambda_{\cdot}|&\leq& 2\frac{1}{N_{1}^{4}}\left( \frac{|v_{\ast\,j}|}{3L}(2\max\{N_{1},4N_{2}\}+3N_{2})\right)^{2}\nonumber\\
&=&\frac{2\,v_{\ast\,j}^{2}}{9L^{2}}\left( \frac{4\max\{N_{1},4N_{2}\}^{2}}{N_{1}^{4}}+6\frac{N_{2}\max\{N_{1},4N_{2}\}}{N_{1}^{4}}+\frac{3N_{2}^{2}}{N_{1}^{4}}\right).\nonumber
\end{eqnarray}
Noticing that this norm will be minimized if $N_{2}<N_{1}$, we conclude the proof.
\end{proof}
\end{proposition}

{\bfseries Remark. }Proposition~\ref{prop:criteria1} only assumes two different representations of the same curve. Nevertheless, the use of the bound found in~(\ref{ec:eigenBound}) assumes that $N_{1}$ and $N_{2}$ do not differ overly because we linearized the inverse of relative distances to get the bounds. Additionally, we suppose that both representations are suitable in the sense that the sum over the polygonal distances approximates the curve length, because the numerical integration depends on this hypothesis in order to handle $v_{\varepsilon}$.

\section{Conclusions}
We presented a novel anisotropic mean curvature geometric flow together with an implementation of it through a numerical scheme, applied to contour recognition tasks. Our proposal has the following features: (1) The curve evolution is given by an explicit scheme. Consequently, the required resolution for the contour recognition can used to choose the size of the stored data array. (2) The stability of our scheme can only be checked at each iteration using the explicit values of $v_{\ast\,j}$ $L$, $N_{1}$ and $N_{2}$, according to Proposition~\ref{prop:criteria1}. (3) We provided specific criteria to improve the conditioning and to verify the stability of this method at each time step.

Moreover, the implementation was optimized to run in parallel, and the code has been made publicly available.

Future improvements for this method may include increasing of the number of point charges that drive the curve evolution, and inferring optimal distributions and values for the potentials. Due to the unavailability of implementations---for example accesible in code repositories---of previous methods, we defer comprehensive comparisons with other schemes of numerical mean curvature flows~\cite{sethian99,sevcovic2011,Mikula}. 

\bibliographystyle{amsplain}

\begin{thebibliography}{10}

\bibitem{Mikula}
Martin Balazovjech and Karol Mikula.
\newblock A higher order scheme for a tangentially stabilized plane curve
  shortening flow with a driving force.
\newblock {\em SIAM Journal on Scientific Computing}, 33(5):2277--2294, 2011.

\bibitem{brakke}
Kenneth~A. Brakke.
\newblock {\em The Motion of a Surface by Its Mean Curvature. (MN-20)}.
\newblock Princeton University Press, 1978.

\bibitem{invpois2}
Martin Burger, Yanina Landa, Nicolay~M. Tanushev, and Richard Tsai.
\newblock {\em Algorithmic Foundation of Robotics VIII: Selected Contributions
  of the Eight International Workshop on the Algorithmic Foundations of
  Robotics}.
\newblock Springer Berlin Heidelberg, Berlin, Heidelberg, 2010.

\bibitem{Chen2006UniquenessAP}
Bing-Long Chen.
\newblock Uniqueness and pseudolocality theorems of the mean curvature flow.
\newblock {\em Communications in Analysis and Geometry}, 15(3):435--490, 2007.

\bibitem{chopp1993}
David~L. Chopp and James~A. Sethian.
\newblock Flow under curvature: singularity formation, minimal surfaces, and
  geodesics.
\newblock {\em Experiment. Math.}, 2(4):235--255, 1993.

\bibitem{chow}
B.~Chow, S.C. Chu, D.~Glickenstein, C.~Guenther, J.~Isenberg, T.~Ivey,
  D.~Knopf, P.~Lu, F.~Luo, and L.~Ni.
\newblock {\em The Ricci Flow: Techniques and Applications:}.
\newblock Mathematical Surveys and Monographs. American Mathematical Society,
  2015.

\bibitem{ecker}
K.~Ecker.
\newblock {\em Regularity Theory for Mean Curvature Flow}.
\newblock Progress in Nonlinear Differential Equations and Their Applications.
  Birkh{\"a}user Boston, 2012.

\bibitem{evans2010partial}
L.C. Evans.
\newblock {\em Partial Differential Equations}.
\newblock Graduate studies in mathematics. American Mathematical Society, 2010.

\bibitem{folland1995introduction}
G.B. Folland.
\newblock {\em Introduction to Partial Differential Equations}.
\newblock Princeton University Press, 1995.

\bibitem{GerhardPolden1999}
Huisken Gerhard and Polden Alexander.
\newblock {\em Calculus of Variations and Geometric Evolution Problems:
  Lectures given at the 2nd Session of the Centro Internazionale Matematico
  Estivo (C.I.M.E.) held in Cetraro, Italy, June 15--22, 1996}.
\newblock Springer Berlin Heidelberg, Berlin, Heidelberg, 1999.

\bibitem{banach}
P.~H{\'a}jek and M.~Johanis.
\newblock {\em Smooth analysis in Banach spaces}.
\newblock De Gruyter Series in Nonlinear Analysis and Applications. De Gruyter,
  2014.

\bibitem{snakes}
Michael Kass, Andrew Witkin, and Demetri Terzopoulos.
\newblock Snakes: Active contour models.
\newblock {\em International Journal of Computer Vision}, 1(4):321--331, 1988.

\bibitem{kimura}
Masato Kimura.
\newblock Numerical analysis of moving boundary problems using the boundary
  tracking method.
\newblock {\em Japan Journal of Industrial and Applied Mathematics},
  14(3):373--398.

\bibitem{krantz}
S.G. Krantz.
\newblock {\em Geometric Function Theory: Explorations in Complex Analysis}.
\newblock Cornerstones. Birkh{\"a}user Boston, 2007.

\bibitem{leveque2007finite}
R.J. LeVeque.
\newblock {\em Finite Difference Methods for Ordinary and Partial Differential
  Equations: Steady-State and Time-Dependent Problems}.
\newblock Society for Industrial and Applied Mathematics, 2007.

\bibitem{invpois3}
Leevan Ling, Y.~C. Hon, and M.~Yamamoto.
\newblock Inverse source identification for poisson equation.
\newblock {\em Inverse Problems in Science and Engineering}, 13(4):433--447,
  2005.

\bibitem{minmaxSethian}
R.~Malladi and J.A. Sethian.
\newblock Image processing: Flows under min/max curvature and mean curvature.
\newblock {\em Graphical Models and Image Processing}, 58(2):127 -- 141, 1996.

\bibitem{mantegazza}
C.~Mantegazza.
\newblock {\em Lecture Notes on Mean Curvature Flow}.
\newblock Progress in Mathematics. Springer Basel, 2011.

\bibitem{met}
W.~W. Mullins.
\newblock Two-dimensional motion of idealized grain boundaries.
\newblock {\em Journal of Applied Physics}, 27(8):900--904, 1956.

\bibitem{invpois1}
Takashi Ohe and Kohzaburo Ohnaka.
\newblock Boundary element approach for an inverse source problem of the
  poisson equation with a one-point-mass like source.
\newblock {\em Applied Mathematical Modelling}, 18(4):216 -- 223, 1994.

\bibitem{osher-level}
S.~Osher and R.~Fedkiw.
\newblock {\em Level Set Methods and Dynamic Implicit Surfaces}.
\newblock Applied Mathematical Sciences. Springer, 2003.

\bibitem{pletcher}
R.H. Pletcher, J.C. Tannehill, and D.~Anderson.
\newblock {\em Computational Fluid Mechanics and Heat Transfer, Third Edition}.
\newblock Series in Computational and Physical Processes in Mechanics and
  Thermal Sciences. CRC Press, 2016.

\bibitem{sauter2010boundary}
S.A. Sauter and C.~Schwab.
\newblock {\em Boundary Element Methods}.
\newblock Springer Series in Computational Mathematics. Springer Berlin
  Heidelberg, 2010.

\bibitem{sethian99}
J.A. Sethian.
\newblock {\em Level Set Methods and Fast Marching Methods: Evolving Interfaces
  in Computational Geometry, Fluid Mechanics, Computer Vision, and Materials
  Science}.
\newblock Cambridge Monographs on Applied and Computational Mathematics.
  Cambridge University Press, 1999.

\bibitem{sevcovic2011}
Daniel {\v{S}}ev{\v{c}}ovi{\v{c}} and Shigetoshi Yazaki.
\newblock Evolution of plane curves with a curvature adjusted tangential
  velocity.
\newblock {\em Japan Journal of Industrial and Applied Mathematics}, 28(3):413,
  2011.

\bibitem{Mcf}
E.~Pedersen T.~Colding, W.~Minicozzi.
\newblock Mean curvature flow.
\newblock {\em Bull. Amer. Math. Soc.}, 52:297 -- 333, 2015.

\bibitem{wu2006elliptic}
Z.~Wu, J.~Yin, and C.~Wang.
\newblock {\em Elliptic \& Parabolic Equations}.
\newblock World Scientific, 2006.

\bibitem{zeidlerFPT}
E.~Zeidler.
\newblock {\em Nonlinear Functional Analysis and Its Applications: Fixed point
  theorems}.
\newblock Nonlinear Functional Analysis and Its Applications. Springer-Verlag,
  1985.

\bibitem{zeidler89}
E.~Zeidler and L.F. Boron.
\newblock {\em Nonlinear Functional Analysis and Its Applications: II/ A:
  Linear Monotone Operators}.
\newblock Monotone operators / transl. by the author and by Leo F. Boron.
  Springer New York, 1989.

\bibitem{zhulectures}
X.P. Zhu.
\newblock {\em Lectures on Mean Curvature Flows}.
\newblock AMS/IP studies in advanced mathematics. American Mathematical Soc.

\end{thebibliography}

\end{document}